\documentclass[a4,11pt]{amsart}
\usepackage{amssymb,amsmath,amsthm,txfonts}
\usepackage{cite}

\newtheorem{thm}{Theorem}[section]
\newtheorem{lem}[thm]{Lemma}

\newtheorem{prop}[thm]{Proposition}
\theoremstyle{definition}

\theoremstyle{remark}
\newtheorem{rem}[thm]{\textbf{Remark}}
\newtheorem{rems}[thm]{\textbf{Remarks}}

 \makeatletter
    
    \@addtoreset{equation}{section}
  \makeatother

\makeatletter
      \def\@makefnmark{%
         \leavevmode
            \raise.9ex\hbox{\check@mathfonts
                \fontsize\sf@size\z@\normalfont%
                            \@thefnmark}%
       }
      \makeatother

\newcommand{\D}{\textrm{div}}

\newcommand{\dd}{\textrm{d}}

\begin{document}

\title[]{Vanishing viscosity limits for axisymmetric flows with boundary}
\author[]{K. Abe}
\date{}
\address[K. ABE]{Department of Mathematics, Graduate School of Science, Osaka City University, 3-3-138 Sugimoto, Sumiyoshi-ku Osaka, 558-8585, Japan}
\email{kabe@sci.osaka-cu.ac.jp}

\date{\today}

\maketitle

\begin{abstract}
We construct global weak solutions of the Euler equations in an infinite cylinder $\Pi=\{x\in \mathbb{R}^{3}\ |\ x_h=(x_1,x_2),\ r=|x_h|<1\}$ for axisymmetric initial data without swirl when initial vorticity $\omega_{0}=\omega^{\theta}_{0}e_{\theta}$ satisfies $\omega^{\theta}_{0}/r\in L^{q}$ for $q\in [3/2,3)$. The solutions constructed are H\"older continuous for spatial variables in $\overline{\Pi}$ if in addition that $\omega^{\theta}_{0}/r\in L^{s}$ for $s\in (3,\infty)$ and unique if $s=\infty$. The proof is by a vanishing viscosity method. We show that the Navier-Stokes equations subject to the Neumann boundary condition is globally well-posed for axisymmetric data without swirl in $L^{p}$ for all $p\in [3,\infty)$. It is also shown that the energy dissipation tends to zero if $\omega^{\theta}_{0}/r\in L^{q}$ for $q\in [3/2,2]$, and Navier-Stokes flows converge to Euler flow in $L^{2}$ locally uniformly for $t\in [0,\infty)$ if additionally $\omega^{\theta}_{0}/r\in L^{\infty}$. The $L^{2}$-convergence in particular implies the energy equality for weak solutions.
\end{abstract}

\vspace{15pt}

\section{Introduction}

\vspace{10pt}
We consider the Navier-Stokes equations:

\begin{equation*}
\begin{aligned}
\partial_t u-\nu\Delta{u}+u\cdot \nabla u+\nabla{p}= 0,\quad \D\ u&=0  \qquad \textrm{in}\ \Pi\times (0,\infty),  \\
\nabla \times u\times n=0,\quad u\cdot n&=0\qquad \textrm{on}\ \partial\Pi\times (0,\infty), \\
u&=u_0\hspace{18pt} \textrm{on}\ \Pi\times\{t=0\},
\end{aligned}
\tag{1.1}
\end{equation*}\\
for the infinite cylinder 

\begin{align*}
\Pi=\{x=(x_1,x_2,x_3)\in \mathbb{R}^{3}\ |\ x_h=(x_1,x_2),\ |x_h|<1  \ \}.
\end{align*}\\
Here, $n$ denotes the unit outward normal vector field on $\partial\Pi$ and $\nu>0$ is the kinematic viscosity.

We study the problem (1.1) for axisymmetric initial data. We say that a vector field  $u$ is axisymmetric if $u(x)={}^{t}Ru(Rx)$ for $x\in \Pi$, $\eta\in [0,2\pi]$, $R=(e_r,e_{\theta},e_z)$ and $e_{r}={}^{t}(\cos\eta,\sin\eta,0)$, $e_{\theta}={}^{t}(-\sin\eta,\cos\eta,0)$, $e_{z}={}^{t}(0,0,1)$. By the cylindrical coordinate $(r,\theta,z)$, an axisymmetric vector field is decomposed into three terms $u=u^{r}e_{r}+u^{\theta}e_{\theta}+u^{z}e_{z}$ and the azimuthal component $u^{\theta}$ is called swirl velocity (e.g., \cite{MaB}). It is known that the Cauchy problem is globally well-posed for axisymmetric initial data without swirl in $H^{2}$ \cite{La68b}, \cite{UI}, \cite{LMNP}. See also \cite{Abidi} for $H^{1/2}$.

The purpose of this paper is to study axisymmetric solutions in $L^{p}$. It is well known that the Cauchy problem is locally well-posed in $L^{p}$ for all $p\geq 3$ \cite{Kato84}. However, global well-posedness results are unknown even if initial data is axisymmetric without swirl. For the two-dimensional case, the problem is globally well-posed in $L^{p}$ for all $p\geq 2$ (including $p=\infty$ \cite{GMS}), since vorticity of two-dimensional flows are uniformly bounded. On the other hand, for axisymmetric flows without swirl, vorticity estimates are more involved due to the vortex stretching as $r\to\infty$.

Recently, global-in-time solutions of the Cauchy problem are constructed in \cite{FengSverak} for axisymmetric data without swirl when initial vorticity $\omega_0=\omega^{\theta}_{0}e_{\theta}$ is a vortex ring, i.e., $\omega^{\theta}_{0}=\kappa\delta_{r_0,z_0}$ for $\kappa\in \mathbb{R}$ and a Dirac measure $\delta_{r_0,z_0}$ in the $(r,z)$-plane. See \cite{GallaySverak2} for the uniqueness. For such initial data, initial velocity belong to $L^{p}$ for $p\in (1,2)$ and $\textrm{BMO}^{-1}$ by the Biot-Savart law. For small data in $\textrm{BMO}^{-1}$, a global well-posedness result is in known \cite{KT01}.

In this paper, we study axisymmetric solutions in the infinite cylinder $\Pi=\{r<1\}$,  subject to the Neumann boundary condition. Since the cylinder is horizontally bounded, vorticity estimates are simpler than those in the whole space. We prove that vorticity of axisymmetric solutions without swirl to (1.1) is uniformly bounded in the infinite cylinder $\Pi$, and unique global-in-time solutions exist for large axisymmetric data without swirl in $L^{p}$ for all $p\in [3,\infty)$.

An important application of our well-poseness result is a vanishing viscosity limit as $\nu \to0$. We apply our global well-posedness result to (1.1) and construct global weak solutions of the Euler equations. Although local well-posedness results are well known for the Euler equations with boundary (see below Theorem 1.1), existence of global weak solutions was unknown. The well-posedness result to (1.1) in $L^{p}$ for $p\in [3,\infty)$ enable us to study weak solutions of the Euler equations when initial vorticity is in $L^{q}$ for $q\in [3/2,3)$ by the Biot-Savart law $1/p=1/q-1/3$.

To state a result, let $L^{p}_{\sigma}$ denote the $L^{p}$-closure of $C_{c,\sigma}^{\infty}$, the space of all smooth solenoidal vector fields with compact support in $\Pi$. The space $L^{p}_{\sigma}$ agrees with the space of all divergence-free vector fields whose normal trace is vanishing on $\partial\Pi$ \cite{ST98}. By a local well-posedness result in the companion paper \cite{A6}, unique local-in-time solutions to (1.1) exist for $u_0\in L^{p}_{\sigma}$ and $p\in [3,\infty)$. Our first result is:

\vspace{10pt}

\begin{thm}
Let $u_0\in L^{p}_{\sigma}$ be an axisymmetric vector field without swirl for $p\in [3,\infty)$. Then, there exists a unique axisymmetric solution without swirl $u\in C([0,\infty); L^{p})\cap C^{\infty}(\overline{\Pi}\times (0,\infty))$ of (1.1) with some associated pressure $p\in C^{\infty}(\overline{\Pi}\times (0,\infty))$. 
\end{thm}

\vspace{15pt}

We apply Theorem 1.1 to construct global weak solutions of the Euler equations:

\begin{equation*}
\begin{aligned}
\partial_t u+u\cdot \nabla u+\nabla{p}= 0,\quad \D\ u&=0  \qquad \textrm{in}\ \Pi\times (0,\infty),  \\
\quad u\cdot n&=0\qquad \textrm{on}\ \partial\Pi\times (0,\infty), \\
u&=u_0\hspace{18pt} \textrm{on}\ \Pi\times\{t=0\}.
\end{aligned}
\tag{1.2}
\end{equation*}

\vspace{10pt}
Unique existence of local-in-time solutions of the Euler equations $u\in C([0,T]; W^{k,q})\cap C^{1}([0,T]; W^{k-1,q})$ is known for sufficiently smooth initial data $u_0\in W^{k,q}$ with integers $k>1+n/q$ and $q\in (1,\infty)$, when $\Pi$ is smoothly bounded in $\mathbb{R}^{n}$ for $n\geq 2$ \cite{EbinMarsden}, \cite{BouBre74}, \cite{Te75}, \cite{KatoLai}. For axisymmetric data without swirl, it is known that local-in-time solutions for $u_0\in H^{s}$ ($s\geq 3$) are continued for all time \cite{Yanagisawa94}. Observe that for $k=2$ and $n=3$, the condition $q\in (3,\infty)$ is required in order to construct local-in-time unique solutions. We construct global weak solutions under the lower regularity condition $q\in [3/2,3)$; see below.

When $\Pi$ is a two-dimensional bounded and simply-connected domain (e.g., a unit disk), global weak solutions of the Euler equations are constructed in \cite{MY92} by a vanishing viscosity method for initial vorticity satisfying $\omega_0\in L^{q}$ and $q\in (1,2)$. For the two-dimensional case, the Neumann boundary condition in (1.1) is reduced to the condition $\omega=0$ and $u\cdot n=0$ on $\partial\Pi$, called the free condition \cite{Lions69} (\cite[p.129]{LionsBook}). The vanishing viscosity method subject to the free condition is studied in \cite{Bardos72}, \cite{VK81} for $\omega_0\in L^{2}$. The condition $q\in (1,2)$ implies that initial velocity belongs to $L^{p}$ for some $p\in (2,\infty)$ by the Biot-Savart law $u_{0}=\nabla^{\perp} (-\Delta_{D})^{-1}\omega_{0}$ for $1/p=1/q-1/2$. Here, $\nabla^{\perp}={}^{t}(\partial_2,-\partial_1)$ and $-\Delta_{D}$ denotes the Laplace operator subject to the Dirichlet boundary condition.

Our goal is to construct three-dimensional weak solutions in the infinite cylinder $\Pi$ for axisymmetric data without swirl when initial vorticity $\omega_{0}=\omega^{\theta}_{0}e_{\theta}$ satisfies $\omega^{\theta}_{0}/r\in L^{q}$ for $q\in [3/2,3)$. The assumption for $\omega^{\theta}_{0}/r$ is stronger than that for vorticity itself and implies that the initial velocity is in $L^{p}$ for some $p\in [3,\infty)$ by the Biot-Savart law $u_0=\nabla \times (-\Delta_{D})^{-1}\omega_{0}$ and $1/p=1/q-1/3$. For such initial data, unique global-in-time solutions to (1.1) exist by  Theorem 1.1. Note that the condition $\omega^{\theta}_{0}/r\in L^{q}$ is weaker than $u_0\in W^{2,q}$ for $q\in [3/2,3)$ since $\omega^{\theta}_{0}/r=-\Delta u^{z}_{0}-\partial_r\omega^{\theta}_{0}$.

Let $BC_{w}([0,\infty); L^{p})$ denote the space of bounded and weakly continuous (resp. weakly-star continuous) functions from $[0,\infty)$ to $L^{p}$ for $p\in (1,\infty)$ (resp. for $p=\infty$). Let $\mathbb{P}$ denote the Helmholtz projection on $L^{p}$ \cite{ST98}. We construct global weak solutions for $\omega^{\theta}_{0}/r\in L^{q}$ and $q\in [3/2,3)$, which are $L^{p}$-integrable and may not be continuous. Under the additional regularity assumptions $\omega^{\theta}_{0}/r\in L^{s}$ for $s\in(3,\infty)$ and $s=\infty$, the weak solutions are H\"older continuous and unique. The main result of this paper is the following:

\vspace{10pt}

\begin{thm}
Let $u_0\in L^{p}_{\sigma}$ be an axisymmetric vector field without swirl for $p\in [3,\infty)$ such that $\omega^{\theta}_{0}/r\in L^{q}$ for $q\in [3/2,3)$ and $1/p=1/q-1/3$.

\noindent
(i) (Existence) There exists a weak solution $u\in BC_{w}([0,\infty); L^{p})$  of (1.2) in the sense that $\nabla u\in BC_{w}([0,\infty); L^{q})$ and

\begin{align*}
\int_{0}^{\infty}\int_{\Pi}(u\cdot\partial_t \varphi +uu:\nabla \varphi)\dd x\dd t
=-\int_{\Pi}u_0\cdot \varphi_{0}\dd x  \tag{1.3}
\end{align*}\\
for all $\varphi\in C^{1}_{c}(\overline{\Pi}\times [0,\infty))$ such that $\D\ \varphi=0$ in $\Pi$ and $\varphi\cdot n=0$ on $\partial\Pi$ for $t\geq 0$, where $\varphi_0(x)=\varphi(x,0)$.

\noindent
(ii) (H\"older continuity) If $\omega^{\theta}_{0}/r\in L^{s}$ for $s\in (3,\infty)$, then  $u\in BC([0,\infty); L^{s})$ satisfies $\nabla u\in BC_{w}([0,\infty); L^{s})$, $\partial_t u\in L^{\infty}(0,\infty; L^{s})$ and 

\begin{align*}
\partial_t u+\mathbb{P}u\cdot \nabla u=0\quad \textrm{on}\ L^{s}\quad \textrm{for a.e.}\ t>0.  \tag{1.4}
\end{align*}\\
In particular, $u(\cdot,t)$ is bounded and H\"older continuous in $\overline{\Pi}$ of exponent $1-3/s$ for each $t\geq 0$. 

\noindent
(iii) (Uniqueness) If in addition that $\omega^{\theta}_{0}/r\in L^{\infty}$,  then $\nabla \times u\in BC_{w}([0,\infty); L^{\infty})$ and the weak solution is unique.
\end{thm}

\vspace{10pt}
It is an interesting question whether the weak solutions constructed in Theorem 1.2 conserve the energy. Since the Poincar\'e inequality holds for the infinite cylinder (see Remarks 4.4 (ii)), the condition $\omega^{\theta}_{0}/r\in L^{q}$ for $q\in [3/2,2]$ implies the finite energy $u_0\in L^{p}_{\sigma}\cap L^{2}$ and the energy equality holds for global-in-time solutions to (1.1); see below (1.5). In the sequel, we consider the case $q\in [3/2,2]$.

In the Kolmogorov's theory of turbulence, it is a basic hypothesis that the energy dissipation tends to a positive constant at large Reynolds numbers. See, e.g., \cite{Frisch}. If the energy dissipation converges to a positive constant for global-in-time solutions $u_{\nu}$ to (1.1) as $\nu\to0$, we would obtain a weak solution strictly decreasing the energy as a vanishing viscosity limit. Unfortunately, due to a regularizing effect, the energy dissipation converges to zero at least under the initial condition $\omega^{\theta}_{0}/r\in L^{q}$ for $q\in [3/2,2]$. See Remarks 2.6 (ii) for $q\in [1,3/2)$. However, it is still non-trivial whether vanishing viscosity limits conserve the energy since they are no longer continuous. 

In the sequel, we prove that global-in-time solutions to (1.1) converge to a limit in $L^{2}$ locally uniformly for $t\in [0,\infty)$ under the additional assumption $\omega^{\theta}_{0}/r\in L^{\infty}$. If the limit is a $C^{1}$-solution, the $L^{2}$-convergence to a limit is equivalent to the convergence of the energy dissipation \cite{Kato84vis}. The equivalence may not always hold if a limit is a weak solution. The assumption $\omega^{\theta}_{0}/r\in L^{\infty}$ does not imply that $\nabla u$ is bounded for the limit and at present is optimal in order to obtain the $L^{2}$-convergence. Once we know the $L^{2}$-convergence, the energy conservation immediately follows as a consequence.

\vspace{15pt}

\begin{thm}
Let $u_{0}\in L^{p}_{\sigma}\cap L^{2}$ be an axisymmetric vector field without swirl such that $\omega_{0}^{\theta}/r\in L^{q}$ for $q\in [3/2,2]$ and $1/p=1/q-1/3$. Let $u_{\nu}$ be a solution of (1.1) in Theorem 1.1.

\noindent 
(i) (Energy dissipation) 
The solution $u_{\nu}\in BC([0,\infty); L^{2})$ satisfies the energy equality 

\begin{align*}
\int_{\Pi}|u_{\nu}|^{2}\dd x+2\nu \int_{0}^{t}\int_{\Pi}|\nabla u_{\nu}|^{2}\dd x\dd s=\int_{\Pi}|u_{0}|^{2}\dd x\qquad t\geq 0,  \tag{1.5}
\end{align*}\\
and 
\begin{align*}
\nu \int_{0}^{T}\int_{\Pi}|\nabla u_{\nu}|^{2}\dd x\dd s=O(\nu^{5/2-3/q})\quad \textrm{as}\ \nu\to0\quad \textrm{for each}\ T>0.  \tag{1.6}
\end{align*}

\noindent
(ii) ($L^{2}$-convergence) Assume in addition that $\omega^{\theta}_{0}/r\in L^{\infty}$. Then, 

\begin{align*}
\lim_{\nu,\ \mu\to0}\sup_{0\leq t\leq T}||u_{\nu}-u_{\mu}||_{L^{2}(\Pi)}=0.   \tag{1.7}
\end{align*}\\
In particular, the limit $u\in BC([0,\infty); L^{2})$ satisfies the energy equality of (1.2):

\begin{align*}
\int_{\Pi}|u|^{2}\dd x=\int_{\Pi}|u_{0}|^{2}\dd x\qquad t\geq 0.  \tag{1.8}
\end{align*}
\end{thm}

\vspace{15pt}
It is noted that there is a possibility that the energy equality (1.8) holds under a weaker assumption than $\omega^{\theta}_{0}/r\in L^{\infty}$ although we assumed it in order to prove the $L^{2}$-convergence (1.7). In fact, it is known as a celebrated Onsager's conjecture \cite{Onsager} that H\"older continuous weak solutions to the Euler equations of exponent $\alpha>1/3$ conserve the energy (but not necessarily if $\alpha\leq 1/3$). The conjecture is studied in \cite{Eyink} and the energy conservation is proved for weak solutions in the whole space under a stronger assumption. A simple proof is given in \cite{CET} under a weaker and natural assumption in the Besov space $u\in L^{3}(0,T; B_{3}^{\alpha,\infty})$ for $\alpha>1/3$. See \cite{DR}, \cite{Constantin08} for further developments and \cite{Eyink06} for a review. Recently, the energy conservation is proved in \cite{BT18} for weak solutions in a bounded domain in the H\"older space $u\in L^{3}(0,T; C^{\alpha}(\overline{\Pi}))$ for $\alpha>1/3$. The weak solutions constructed in Theorem 1.2 are indeed H\"older continuous of exponent $\alpha=1-3/s>1/3$ if in addition that $\omega^{\theta}_{0}/r\in L^{s}$ for $s\in (9/2,\infty]$. If the result of \cite{BT18} holds also for the infinite cylinder, the weak solutions in Theorem 1.2 satisfy (1.8) even for $s\in (9/2,\infty]$.

\vspace{15pt}

We outline the proofs of Theorems 1.1-1.3. By a local well-posedness result of (1.1) in \cite{A6}, there exist local-in-time smooth axisymmetric solutions without swirl $u\in C([0,T]; L^{p})\cap C^{\infty}(\overline{\Pi}\times (0,T])$ for $u _0\in L^{p}_{\sigma}$ and $p\in [3,\infty)$ satisfying the integral equation

\begin{align*}
u=e^{-t\nu A}u_0-\int_{0}^{t}e^{-(t-s)\nu A}\mathbb{P}\ (u\cdot \nabla u)(s)\dd s.   
\end{align*}\\
Here, $A$ denotes the Stokes operator subject to the Neumann boundary condition. We establish an apriori estimate in $L^{p}$ based on the vorticity equation. Since the vorticity $\omega=\omega^{\theta}e_{\theta}$ vanishes on the boundary subject to the Neumann boundary condition, $\omega^{\theta}/r$ satisfies the drift-diffusion equation with the homogeneous Dirichlet boundary condition:

\begin{equation*}
\begin{aligned}
\partial_t \Big(\frac{\omega^{\theta}}{r}\Big)+u\cdot \nabla\Big(\frac{\omega^{\theta}}{r}\Big)-\nu \Big(\Delta+\frac{2}{r}\partial_r\Big)\Big(\frac{\omega^{\theta}}{r}\Big)&=0\quad \textrm{in}\ \Pi\times (0,T),\\
\frac{\omega^{\theta}}{r}&=0\quad \textrm{on}\ \partial\Pi\times (0,T).
\end{aligned}
\tag{1.9}
\end{equation*}\\
We prove the a priori estimate

\begin{align*}
\Big\|\frac{\omega^{\theta}}{r}\Big\|_{L^{r}(\Pi)}\leq \frac{C}{(\nu t)^{\frac{3}{2}(\frac{1}{q}-\frac{1}{r}) }}\Big\|\frac{\omega^{\theta}_{0}}{r}\Big\|_{L^{q}(\Pi)}\qquad t>0,\ \nu>0,  \tag{1.10}
\end{align*}\\
for $1\leq q\leq r\leq \infty$. The estimate (1.10) for $r=q$ is proved in \cite{UI} for $\Pi=\mathbb{R}^{3}$. Moreover, the decay estimate for $r\in [1,\infty]$ and $q=1$ is established in \cite{FengSverak}. We prove (1.10) for the infinite cylinder $\Pi$. Since local-in-time solutions of (1.1) are smooth for $t>0$, we may assume that $\omega^{\theta}_{0}/r$ is bounded. Then, the a priori estimate (1.10) for $r=q=\infty$ implies that vorticity is uniformly bounded in the infinite cylinder $\Pi=\{r<1\}$. Since $\mathbb{P}u\cdot \nabla u=\mathbb{P}\omega\times u$, by the Gronwall's inequality we obtain an exponential bound of the form

\begin{align*}
||u||_{L^{p}(\Pi)}\leq C||u_0||_{L^{p}(\Pi)} \exp\Big(C\Big\|\frac{\omega^{\theta}_0}{r}\Big\|_{L^{\infty}(\Pi)} t\Big)\quad t\geq 0.  
\end{align*}\\
The estimate implies that the $L^{p}$-norm does not blow-up. Hence the local-in-time solutions are continued for all time. (Moreover, the solutions converge to zero in $L^{r}$ for $r\in (p,\infty)$ as time goes to infinity; see Remarks 2.6 (iii).)

The proof of Theorem 1.2 is based on the Biot-Savart law in the infinite cylinder. We show that axisymmetric vector fields without swirl satisfy

\begin{align*}
u=\nabla \times (-\Delta_{D})^{-1} (\nabla\times u),   \tag{1.11}
\end{align*}
\begin{align*}
||u||_{L^{p}}+||\nabla u||_{L^{q}}\leq C||\nabla \times u||_{L^{q}}, \quad 1/p=1/q-1/3.    \tag{1.12}
\end{align*}\\
The existence of global weak solutions (i) and regularity properties (ii) follow from the a priori estimate (1.10) for $r=q$ and (1.12) by taking a vanishing viscosity limit and applying an abstract compactness theorem. The uniqueness in Theorem 1.2 (iii) is based on the growth estimate of the $L^{r}$-norm

\begin{align*}
||\nabla u||_{L^{r}(\Pi)}&\leq Cr||\nabla \times u||_{L^{r}\cap L^{r_0}(\Pi)},  \tag{1.13} 
\end{align*}\\
for $3<r_0< r<\infty$ with some absolute constant $C$. The estimate (1.13) is proved in \cite{Yudovich62} for bounded domains. We extend it for the infinite cylinder and adjust the Yudovich's energy method of uniqueness \cite{Yudovich63} for solutions with infinite energy by a cut-off function argument.

The convergence of the energy dissipation (1.6) follows from the vorticity estimate (1.10) for $r=2$ and $q\in [3/2,2]$. The $L^{2}$-convergence (1.7) is based on the estimate (1.13). Since the condition $\omega^{\theta}_{0}/r\in L^{\infty}$ implies that $\nabla \times u_{\nu}$ is uniformly bounded for $\nu>0$ and $r>3$, we estimate the energy norm of $u_{\nu}-u_{\mu}$ for two solutions of (1.1) by using the estimate (1.13).  

\vspace{20pt}

This paper is organized as follows. In Section 2, we prove the vorticity estimate (1.10) for local-in-time solutions to (1.1). Since we only use the estimate (1.10) for $r=q=\infty$ in order to prove Theorem 1.1, we give a proof for the case $r\neq q$ in Appendix A. In Section 3, we prove the Biot-Savart law (1.11) and the estimate (1.13). In Section 4, we prove Theorem 1.2 (i) and (ii) by applying a vanishing viscosity method. In Section 5, we prove Theorem 1.2 (iii). In Section 6, we prove Theorem 1.3.

\vspace{20pt}

\section{Global smooth solutions with viscosity}

\vspace{15pt}
We prove Theorem 1.1. We first observe unique existence of local-in-time axisymmetric solutions without swirl for $u_0\in L^{p}_{\sigma}$ and $p\in [3,\infty)$. The vorticity estimate (1.10) for $r=q$ is obtained by integration by parts for $q\in [1,\infty)$ and a maximum principle for $q=\infty$. Throughout this Section, we denote solutions of (1.1) by $u=u_{\nu}$ and suppressing $\nu>0$.

\subsection{Local-in-time solutions}

\vspace{15pt}

We set the Laplace operator subject to the Neumann boundary condition  

\begin{align*}
&Bu=-\Delta u,\quad \textrm{for}\ u\in D(B),\\
&D(B)=\{u\in W^{2,p}(\Pi)\ |\ \nabla\times \ u\times n=0,\ u\cdot n=0\ \textrm{on}\ \partial\Pi\  \},
\end{align*}\\
It is proved in \cite[Lemma B.1]{A6} that the operator $-B$ generates a bounded $C_0$-analytic semigroup on $L^{p}$ ($1<p<\infty$) for the infinite cylinder $\Pi$. We set the Stokes operator 

\begin{align*}
&Au=B u,\quad \textrm{for}\ u\in D(A),\\
&D(A)=L^{p}_{\sigma}\cap D(B).
\end{align*}\\
Since $Au\in L^{p}_{\sigma}$ by the Neumann boundary condition, the operator $-A$ generates a bounded $C_0$-analytic semigroup on the solenoidal vector space $L^{p}_{\sigma}$. By the analyticity of the semigroup, we are able to construct local-in-time solutions satisfying the integral form

\begin{align*}
u=e^{-t\nu A}u_0-\int_{0}^{t}e^{-(t-s)\nu A}\mathbb{P}\ (u\cdot \nabla u)(s)\dd s.   \tag{2.1}
\end{align*}\\
By a standard argument using a fractional power of the Stokes operator, it is not difficult to see that all derivatives of the mild solution belong to the H\"older space $C^{\mu}((0,T]; L^{s})$ for $\mu\in (0,1/2)$ and $s\in (3,\infty)$. Hence the mild solution is smooth for $t>0$ and satisfies (1.1).

\vspace{15pt}

\begin{lem}
For an axisymmetric vector field without swirl $u_0\in L^{p}_{\sigma}$ and $p\in [3,\infty)$, there exists $T>0$ and a unique axisymmetric mild solution without swirl $u\in C([0,T]; L^{p})\cap C^{\infty}(\overline{\Pi}\times (0,T])$ of (1.1).
\end{lem}

\vspace{5pt}

\begin{proof}
The unique existence of local-in-time smooth mild solutions is proved in \cite[Theorem 1.1]{A6}. The axial symmetry follows from the  uniqueness. We consider a rotation operator $U: f\longmapsto {}^{t}Rf(Rx)$ for $R=(e_{r}(\eta), e_{\theta}(\eta),e_{z})$ and $\eta\in [0,2\pi]$. Since the Stokes semigroup $e^{-t\nu A}$ and the Helmholtz projection $\mathbb{P}$ are commutable with the operator $U$ (see \cite[Proposition 2.6]{AS}), by multiplying $U$ by (2.1) we see that $Uu$ is a mild solution for the same axisymmetric initial data $u_0$. By the uniqueness of mild solutions, the function $Uu$ agrees with $u$. Hence $u(x,t)={}^{t}Ru(Rx,t)$ for $\eta\in [0,2\pi]$ and $u$ is axisymmetric.

It is not difficult to see that $u$ is without swirl. By the Neumann boundary condition in (1.1), we see that an axisymmetric solution $u=u^{r}e_{r}+u^{\theta}e_{\theta}+u^{z}e_{z}$ satisfies 

\begin{align*}
u^{r}=0,\quad \partial_r u^{\theta}+u^{\theta}=0,\quad \partial_r u^{z}=0\qquad \textrm{on}\ \{r=1\}.    \tag{2.2}
\end{align*}\\
Since $u$ is smooth for $t>0$ and $u^{\theta}_{0}=0$, by a fundamental calculation, we see that $\varphi=u^{\theta}e_{\theta}\in C([0,T]; L^{p})\cap C^{\infty}(\overline{\Pi}\times (0,T))$  satisfies 

\begin{align*}
\partial_t \varphi-\Delta \varphi +u\cdot \nabla \varphi-(\partial_r u^{r}+ \partial_z u^{z})\varphi&=0\quad \textrm{in}\ \Pi\times (0,T),\\
\partial_n\varphi+\varphi&=0\quad \textrm{on}\ \partial\Pi\times (0,T),\\
\varphi&=0\quad \textrm{on}\ \Pi\times \{t=0\}.
\end{align*}\\
Since the Laplace operator $-\Delta_{R}$ with the Robin boundary condition generates a $C_0$-analytic semigroup on $L^{p}$ \cite{ADN} (\cite[Theorem 3.1.3]{Lunardi}), by the uniqueness of the inhomogeneous heat equation, the function $\varphi$ satisfies the integral form 

\begin{align*}
\varphi=-\int_{0}^{t}e^{(t-s)\nu \Delta_{R}}(u\cdot \nabla \varphi-(\partial_r u^{r}+ \partial_z u^{z})\varphi)\dd s.
\end{align*}\\
Since $u\in C([0,T]; L^{p})$ and $t^{1/2}\nabla u\in C([0,T]; L^{p})$, it is not difficult to show that $\varphi\equiv 0$ by estimating $L^{p}$-norms of $\varphi$. Hence the local-in-time solution $u$ is axisymmetric without swirl. 
\end{proof}

\vspace{15pt}
In order to prove Theorem 1.3 later in Section 6, we show that local-in-time solutions satisfy the energy equality (1.5) for initial data with finite energy. 
\vspace{15pt}

\begin{prop}
For axisymmetric initial data without swirl $u_0\in L^{p}_{\sigma}\cap L^{2}$ for $p\in [3,\infty)$, the local-in-time solution $u$ satisfies

\begin{align*}
u,\ t^{1/2}\nabla u\in C([0,T]; L^{p}\cap L^{2}),    \tag{2.3}
\end{align*}\\
and the energy equality (1.5) for $t\geq 0$.
\end{prop}

\vspace{5pt}

\begin{proof}
We prove (2.3). The energy equality (1.5) follows from (2.2) and integration by parts. We give a proof for the case $p=3$ since we are able to prove the case $p\in (3,\infty)$ by a similar way. We may assume that $\nu =1$. We invoke an iterative argument in \cite[Theorem 5.2]{A6}. We use regularizing estimates of the Stokes semigroup \cite[Lemma 5.1]{A6},

\begin{align*}
||\partial_x^{k}e^{-tA}f||_{L^{2}}\leq \frac{C}{t^{\frac{3}{2}(\frac{1}{r}-\frac{1}{2})+\frac{|k|}{2}  }}||f||_{L^{r}}   \tag{2.4}
\end{align*}\\
for $t\leq T_0,\ |k|\leq 1$ and $r\in [6/5,2]$. We set a sequence $\{u_j\}$ as usual by $u_1=e^{-tA}u_0$,

\begin{align*}
u_{j+1}=e^{-tA}u_0-\int_{0}^{t}e^{-(t-s)A}\mathbb{P}u_{j}\cdot \nabla u_{j}\dd s,\quad j\geq 1,
\end{align*} \\
and the constants

\begin{align*}
K_{j}&=\sup_{0\leq t\leq T}t^{\gamma} ||u_{j}||_{L^{q}}(t) ,\\
M_{j}&=\sup_{0\leq t\leq T}( ||u_{j}||_{L^{2}}+t^{1/2}||\nabla u_{j}||_{L^{2}} ),
\end{align*}\\
for $\gamma =3/2(1/3-1/q)$ and $q\in (3,\infty)$. Then, we have $K_j\leq K_1$ for all $j\geq 1$ for sufficiently small $T>0$. We set $1/r=1/2+1/q$. Since $r\in (6/5,2]$, applying (2.4) and the H\"older inequality imply that 

\begin{align*}
||u_{j+1}||_{L^{2}}
&\leq ||e^{-tA}u_0||_{L^{2}}+\int_{0}^{t}\frac{C}{(t-s)^{\frac{3}{2}(\frac{1}{r}-\frac{1}{2}) }}||u_j\cdot \nabla u_j||_{L^{r}}\dd s\\
&\leq ||e^{-tA}u_0||_{L^{2}}+C'K_jM_j.
\end{align*}\\
We estimate the $L^{2}$-norm of $\nabla u_{j+1}$ in a similar way and obtain 

\begin{align*}
M_{j+1}\leq M_1+CK_1M_j.
\end{align*}\\
We take $T>0$ sufficiently small so that $CK_1\leq 1/2$ and obtain the uniform bound $M_{j+1}\leq 2M_1$ for all $j\geq 1$. Since the sequence $\{u_j\}$ converges to a limit $u\in C([0,T]; L^{3})$ such that $t^{1/2}\nabla u \in C([0,T]; L^{3})$, in a similar way, the uniform estimate for $M_j$ is inherited by the limit. We obtained (2.3). 
\end{proof}

\vspace{15pt}

\subsection{Vorticity estimates}

We shall prove the vorticity estimate (1.10) for $r=q\in [1,\infty]$. We first show the case $q\in [1,\infty)$ by integration by parts.

\vspace{15pt}

\begin{lem}
Let $u$ be an axisymmetric solution in Lemma 2.1. Assume that $\omega^{\theta}_{0}/r\in L^{q}$ for $q\in [1,\infty]$. Then, the estimate 

\begin{align*}
\Big\|\frac{\omega^{\theta}}{r}\Big\|_{L^{q}(\Pi)}\leq \Big\|\frac{\omega^{\theta}_{0}}{r}\Big\|_{L^{q}(\Pi)}\quad t>0  \tag{2.5}
\end{align*}\\
holds.
\end{lem}

\vspace{15pt}

\begin{prop}
The estimate (2.5) holds for $q\in [1,\infty)$.
\end{prop}

\vspace{5pt}

\begin{proof}
We observe that $\omega^{\theta}/r$ is smooth for $t>0$ and satisfies 

\begin{equation*}
\begin{aligned}
\partial_t \Big(\frac{\omega^{\theta}}{r}\Big)+u\cdot \nabla\Big(\frac{\omega^{\theta}}{r}\Big)-\nu \Big(\Delta+\frac{2}{r}\partial_r\Big)\Big(\frac{\omega^{\theta}}{r}\Big)&=0\quad \textrm{in}\ \Pi\times (0,T),\\
\Big(\frac{\omega^{\theta}}{r}\Big)&=0\quad \textrm{on}\ \partial\Pi\times (0,T).
\end{aligned}
\tag{2.6}
\end{equation*}\\
In order to differentiate the $L^{q}$-norm of $\Omega=\omega^{\theta}/r$, we approximate the absolute value function $\psi(s)=|s|$. For an arbitrary $\varepsilon>0$, we set a smooth non-negative convex function $\psi_{\varepsilon}(s)=(s^{2}+\varepsilon^{2})^{1/2}-\varepsilon$ for $s\in \mathbb{R}$, i.e., $0\leq \psi_{\varepsilon}\leq |s|$, $\ddot{\psi}_{\varepsilon}>0$. The function $\psi_{\varepsilon}$ satifies $\psi_{\varepsilon}(0)=\dot{\psi}_{\varepsilon}(0)=0$. We differentiate $\psi_{\varepsilon}^{q}(\Omega)$ to see that 

\begin{align*}
\partial_t \psi^{q}_{\varepsilon}(\Omega)&=q \partial_t \Omega \dot{\psi}_{\varepsilon}(\Omega){\psi}^{q-1}_{\varepsilon}(\Omega),\\
\nabla \psi^{q}_{\varepsilon}(\Omega)&=q \nabla \Omega \dot{\psi}_{\varepsilon}(\Omega){\psi}^{q-1}_{\varepsilon}(\Omega).
\end{align*}\\
Since $\psi_{\varepsilon}(\Omega)=\dot{\psi}_{\varepsilon}(\Omega)=0$ on $\partial\Pi$ by the boundary condition, integration by parts yields

\begin{align*}
\frac{d}{dt}\int_{\Pi}\psi_{\varepsilon}^{q}(\Omega)\dd x
&=q\int_{\Pi}\big(-u\cdot \nabla \Omega+\nu \Delta \Omega+2\nu \frac{1}{r}\partial_r\Omega\big)\dot{\psi}_{\varepsilon}(\Omega)\psi^{q-1}_{\varepsilon}(\Omega)\dd x\\
&=-\int_{\Pi}u\cdot \nabla \psi_{\varepsilon}^{q}(\Omega)\dd x+\nu q\int_{\Pi} \Delta \Omega\dot{\psi}_{\varepsilon}(\Omega)\psi^{q-1}_{\varepsilon}(\Omega)\dd x+2\nu\int_{\Pi} \frac{1}{r}\partial_r\psi^{q}_{\varepsilon}(\Omega)\dd x.
\end{align*}\\
The first-term vanishes by the divergence-free condition. Since $\psi_{\varepsilon}$ is non-negative and convex, we see that 

\begin{align*}
\int_{\Pi} \Delta \Omega\dot{\psi}_{\varepsilon}(\Omega)\psi^{q-1}_{\varepsilon}(\Omega)\dd x
&=-\int_{\Pi}|\nabla \Omega|^{2}\big((q-1)\psi^{q-2}_{\varepsilon}(\Omega)|\dot{\psi}_{\varepsilon}(\Omega)|^{2}+\psi_{\varepsilon}^{q-1}(\Omega)\ddot{\psi}_{\varepsilon}(\Omega)  \big)\dd x\\
&\leq -\frac{4}{q}\Big(1-\frac{1}{q}\Big)\int_{\Pi}\big|\nabla \psi_{\varepsilon}(\Omega)^{\frac{q}{2}}\big|^{2}\dd x,
\end{align*}
\begin{align*}
&\int_{\Pi}\frac{1}{r}\partial_r\psi^{q}_{\varepsilon}(\Omega)\dd x
=2\pi \int_{\mathbb{R}}\dd z\int_{0}^{1}\partial_r\psi^{q}_{\varepsilon}(\Omega)\dd r 
=-2\pi \int_{\mathbb{R}}\psi^{q}_{\varepsilon}(\Omega(0,z,t))\dd z\leq 0,
\end{align*}\\
for $\Omega=\Omega(r,z,t)$. Hence we have

\begin{align*}
\frac{d}{dt}\int_{\Pi}\psi_{\varepsilon}^{q}(\Omega)\dd x
+4\nu\Big(1-\frac{1}{q}\Big)\int_{\Pi}\big|\nabla \psi_{\varepsilon}(\Omega)^{\frac{q}{2}}\big|^{2}\dd x
\leq 0.    \tag{2.7}
\end{align*}\\
We integrate in $[0,t]$ and estimate

\begin{align*}
\int_{\Pi}\psi_{\varepsilon}^{q}(\Omega)\dd x\leq \int_{\Pi}\psi_{\varepsilon}^{q}(\Omega_0)\dd x \quad t\geq 0.
\end{align*}\\
Since $\psi_{\varepsilon}(s)$ monotonically converges to $\psi(s)=|s|$, sending $\varepsilon\to 0$ implies the desired estimate for $q\in [1,\infty)$.
\end{proof}

\vspace{15pt}
Following \cite{KNSS}, \cite{FengSverak}, we prove the case $q=\infty$ by a maximum principle.

\vspace{15pt}

\begin{proof}[Proof of Lemma 2.3]
We apply a maximum principle for $\Omega=\omega^{\theta}/r$ by regarding $\Delta_{x}+2r^{-1}\partial_r$ as the Laplace operator in $\mathbb{R}^{5}$. Let $S^{3}$ denote the unit sphere in $\mathbb{R}^{4}$. Let $(r,\theta,z)$ be the cylindrical coordinate for the cartesian coordinate $x=(x_1,x_2,x_3)$. For $r>0$, $z\in \mathbb {R}$ and $\tau \in S^{3}$, we set  new variables $y=(y_{h},y_5)$ and $y_h=(y_1,y_2,y_3,y_4)$ by $y_h=r\tau$ and $y_5=z$. Then, the gradient and the Laplace operator are written as $\nabla_{y}=\tau \partial_r +\nabla_{S^{3}}+e_{5}\partial_z$ and $\Delta_{y}=\partial_r^{2}+3r^{-1}\partial_r+\Delta_{S^{3}}+\partial_{z}^{2}$ with the surface gradient $\nabla_{S^{3}}$. Let $e_5={}^{t}(0,0,0,0,1)$. We define $\tilde{u}(y,t)$ and $\tilde{\Omega}(y,t)$ by

\begin{align*}
&\tilde{u}(y,t)=u^{r}(r,z,t)\tau+u^{z}(r,z,t)e_{5},\\
&\tilde{\Omega}(y,t)=\Omega(r,z,t).
\end{align*}\\
Since $\tilde{u}\cdot \nabla_y=u^{r}\partial_r+u^{z}\partial_z=u\cdot \nabla_x$ and $\Delta_{y}=\Delta_{x}+2r^{-1}\partial_r+\Delta_{S^{3}}$, the vorticity equation (2.6) is then written as 

\begin{equation*}
\begin{aligned}
\partial_t \tilde{\Omega}+\tilde{u}\cdot \nabla_y \tilde{\Omega}-\nu\Delta_{y}\tilde{\Omega}&=0\quad \textrm{in}\ \Pi_5\times (0,T),\\
\tilde{\Omega}&=0\quad \textrm{on}\ \partial\Pi_5\times (0,T),
\end{aligned}
\tag{2.8}
\end{equation*}\\
for the five-dimensional cylinder $\Pi_5=\{y\in \mathbb{R}^{5}\ |\ |y_h|<1\}$. 

Since $\tilde{\Omega}$ may not be continuous at $t=0$, we approximate initial data by $u_{0,\varepsilon}=e^{-\varepsilon A}u_0$, $\varepsilon>0$, and apply a maximum principle for solutions to $u_{0,\varepsilon}$. Let $-\Delta_{D}$ denote the Laplace operator subject to the Dirichlet boundary condition in the cylinder $\Pi_{5}$. Since $u_0$ is axisymmetric without swirl and $\omega^{\theta}_{0}/r\in L^{\infty}$, by the uniqueness of the heat equation, we see that $\omega^{\theta}_{0,\varepsilon}/r=e^{\varepsilon \Delta_{D}}(\omega^{\theta}_{0}/r)$. Since the heat semigroup is a contraction semigroup on $L^{\infty}$, it follows that 

\begin{align*}
\Big\|\frac{\omega^{\theta}_{0,\varepsilon}}{r}\Big\|_{L^{\infty}(\Pi)}\leq
\Big\|\frac{\omega^{\theta}_{0}}{r}\Big\|_{L^{\infty}(\Pi)}.
\end{align*}\\
By Lemma 2.1, there exists a local-in-time axisymmetric solution without swirl $u_{\varepsilon}$ of (1.1) for $u_{0,\varepsilon}$. Since $u_{0,\varepsilon}\in D(A^m_{s})$ for $s\in (3,\infty)$ and all $m\geq 0$, all derivatives of the local-in-time solution are bounded and continuous up to time zero \cite[Remarks 6.5(i)]{A6}. Here, $A_{s}$ denotes the Stokes operator in $L^{s}_{\sigma}$. Hence applying the maximum principle for $(\tilde{u}_{\varepsilon}, \tilde{\Omega}_{\varepsilon})$ yields 

\begin{align*}
\Big\|\frac{\omega^{\theta}_{\varepsilon}}{r}\Big\|_{L^{\infty}(\Pi)}\leq
\Big\|\frac{\omega^{\theta}_{0,\varepsilon}}{r}\Big\|_{L^{\infty}(\Pi)}
\leq \Big\|\frac{\omega^{\theta}_{0}}{r}\Big\|_{L^{\infty}(\Pi)}.
\end{align*}\\
Since $u_{0,\varepsilon}$ converges to $u_0$ in $L^{p}$, the local-in-time  solution $u_{\varepsilon}$ converges to the mild solution $u\in C([0,T]; L^{p})$ for $u_0$ and the vorticity estimate is inherited to the limit. Thus (2.5) holds for $q=\infty$.
\end{proof}

\vspace{15pt}

We further deduce the decay estimate of vorticity from the inequality (2.7). We give a proof for the following Lemma 2.5 in Appendix A.   

\vspace{15pt}

\begin{lem}
Under the same assumption of Lemma 2.3, the estimate 

\begin{align*}
\Big\|\frac{\omega^{\theta}}{r}\Big\|_{L^{r}(\Pi)}\leq \frac{C}{(\nu t)^{\frac{3}{2}(\frac{1}{q}-\frac{1}{r}) }}\Big\|\frac{\omega^{\theta}_{0}}{r}\Big\|_{L^{q}(\Pi)}\qquad t>0,\ \nu>0,  \tag{2.9}
\end{align*}\\
holds for $1\leq q\leq r\leq \infty$ with some constant $C$.
\end{lem}

\vspace{15pt}

\subsection{An exponential bound}

We now complete: 
\vspace{15pt}

\begin{proof}[Proof of Theorem 1.1]
Let $u\in C([0,T]; L^{p})\cap C^{\infty}(\overline{\Pi}\times (0,T])$ be a local-in-time axisymmetric solution in Lemma 2.1. By replacing the initial time to some $t_0\in (0,T]$, we may assume that $\omega^{\theta}_{0}/r\in L^{\infty}$. We apply Lemma 2.3 and estimate the $L^{\infty}$-estimate of $\omega^{\theta}/r$. Since $r<1$, we have

\begin{align*}
\|\omega^{\theta}\|_{L^{\infty}(\Pi)}\leq \Big\|\frac{\omega^{\theta}_{0}}{r}\Big\|_{L^{\infty}(\Pi)}\quad t\geq 0.
\end{align*}\\
Since $\mathbb{P}u\cdot \nabla u=\mathbb{P}\omega\times u$ and the Stokes semigroup is a bounded semigroup on $L^{p}$, it follows from (2.1) that

\begin{align*}
||u||_{L^{p}(\Pi)}
\leq C_1||u_0||_{L^{p}(\Pi)}+C_2\Big\|\frac{\omega^{\theta}_{0}}{r}\Big\|_{L^{\infty}(\Pi)}\int_{0}^{t}||u||_{L^{p}(\Pi)}\dd s\qquad t\geq0.
\end{align*}\\
with some constants $C_1$ and $C_2$, independent of the viscosity $\nu$. Applying the Gronwall's inequality yields 

\begin{align*}
||u||_{L^{p}(\Pi)}\leq C_1||u_0||_{L^{p}(\Pi)} \exp\Big(C_2\Big\|\frac{\omega^{\theta}_0}{r}\Big\|_{L^{\infty}(\Pi)} t\Big)\quad t\geq 0.  \tag{2.10}
\end{align*}\\
Since the $L^{p}$-norm of $u$ is globally bounded, the local-in-time solution is continued for all $t>0$. The proof is complete.
\end{proof}

\vspace{15pt}

\begin{rems}
\noindent 
(i) ($p=\infty$)
It is unknown whether the assertion of Theorem 1.1 holds for $p=\infty$. For the two-dimensional Cauchy problem, unique existence of global-in-time solutions is known for bounded and non-decaying initial data $u_0\in L^{\infty}_{\sigma}$ \cite{GMS}. Moreover, global-in-time solutions satisfy a single exponential bound of the form 

\begin{align*}
||u||_{L^{\infty}(\mathbb{R}^{2})}\leq C_1||u_0||_{L^{\infty}(\mathbb{R}^{2})} \exp\big(C_2||\omega_0||_{L^{\infty}(\mathbb{R}^{2})}t\big)\qquad t\geq 0,
\end{align*}\\
with some constants $C_1$ and $C_2$, independent of viscosity \cite{ST07}. The single exponential bound is further improved to a linear growth estimate as $t\to\infty$ by using viscosity. See \cite{Zelik13}, \cite{GallayLec}. Note that for a two-dimensional layer, unique global-in-time solutions exist for bounded initial data subject to the Neumann boundary condition. Moreover, the $L^{\infty}$-norm of solutions are uniformly bounded for all time \cite{GallaySl}, \cite{GallaySl2}.

\noindent 
(ii) ($1\leq q<3/2$)
It is unknown whether unique global-in-time solutions to (1.1) exist for $\omega^{\theta}_{0}/r\in L^{q}$ and $q\in [1,3/2)$. This condition implies that initial velocity belongs to $L^{p}$ for $p\in [3/2,3)$ and $1/p=1/q-1/3$ by the Biot-Savart law. Although Theorem 1.1 may not be available for this case, it is still likely that unique global-in-time solutions to (1.1) exist. In fact, Gallay-{\v{S}}ver\'ak \cite{GallaySverak} constructed unique global-in-time solutions of the Cauchy problem for $\omega^{\theta}_{0}/r\in L^{1}$, based on the a priori estimate (2.5) for $r=q=1$ and the vorticity equation in a half plane. Note that the convergence of the energy dissipation (1.6) does not follow from the vorticity estimate (2.9) if $q\in [1,6/5]$.

\noindent 
(iii) (Large time behavior)
Global-in-time solutions in Theorem 1.1 are uniformly bounded for all time, i.e., $u\in BC([0,\infty); L^{p})$. In fact, we are able to assume that $\omega^{\theta}_{0}/r\in L^{p}$ by replacing the initial time since local-in-time solutions $u(\cdot,t)$ belong to $W^{2,p}$ for $p\in [3,\infty)$. As proved later in Lemma 3.5 and Proposition 4.3, since axisymmetric solutions of (1.1) are uniquely determined by the Biot-Savart law and the Poincar\'e inequality holds in the cylinder, we have

\begin{align*}
||u||_{L^{p}(\Pi)}
\leq C||\nabla u||_{L^{p}(\Pi)}
\leq C' ||\nabla \times u||_{L^{p}(\Pi)}
\leq C'\Big\| \frac{\omega^{\theta}}{r}\Big\|_{L^{p}(\Pi)}.
\end{align*}\\
By the vorticity estimate (2.9), the solutions are uniformly bounded in $L^{p}$ and tend to zero in $L^{r}$ for $r \in (p,\infty)$ as $t\to\infty$.
\end{rems}

\vspace{15pt}

\section{The Biot-Savart law}

\vspace{15pt}

In this section, we give a Biot-Savart law in the infinite cylinder (Lemma 3.5). Since stream functions exist for axisymmetric vector fields without swirl and satisfy the Dirichlet boundary condition, we are able to represent axisymmetric vector fields without swirl by the Laplace operator with the Dirichlet boundary condition $u=\nabla \times (-\Delta_{D})^{-1}\nabla \times u$. We first prepare $L^{p}$-estimtates for the Dirichlet problem of the Poisson equation and apply them to axisymmetric vector fields without swirl.

\vspace{15pt}

\subsection{$L^{p}$-estimates for the Poisson equation}

We consider the Poisson equation in the infinite cylinder:

\vspace{15pt}

\begin{equation*}
\begin{aligned}
-\Delta \phi=f\quad \textrm{in}\ \Pi,\qquad \phi=0\quad \textrm{on}\ \partial\Pi.
\end{aligned}
\tag{3.1}
\end{equation*}

\vspace{15pt}

\begin{lem}
\noindent
(i) Let $q\in (1,\infty)$. For $f\in L^{q}$, there exists a unique solution $\phi\in W^{2,q}$ of (3.1) satisfying 

\begin{align*}
||\phi||_{W^{2,q}}\leq C||f||_{L^{q}}    \tag{3.2}
\end{align*}\\
with some constant $C$.

\noindent
(ii) For $q\in (1,3)$ and $p\in (3/2,\infty)$ satisfying $1/p=1/q-1/3$, there exists a constant $C'$ such that 

\begin{align*}
||\nabla \phi||_{L^{p}}\leq C'||f||_{L^{q}}.   \tag{3.3}
\end{align*}
\end{lem}

\vspace{15pt}
We prove Lemma 3.1 by using the heat semigroup $e^{t\Delta_{D}}$. 
\vspace{15pt}

\begin{prop}
\noindent 
(i) There exists a constant $M$ such that 

\begin{align*}
||e^{t\Delta_{D}}f||_{L^{q}}\leq e^{-\mu_q t}||f||_{L^{q}}\quad t>0,  \tag{3.4}
\end{align*}\\
for $f\in L^{q}$ with the constant $\mu_q=M/qq'$, where $q'$ is the  conjugate exponent to $q\in (1,\infty)$.

\noindent 
(ii) The heat kernel $K(x,y,t)$ of $e^{t\Delta_{D}}$ satisfies the Gaussian upper bound,

\begin{align*}
0\leq K(x,y,t)\leq \frac{1}{(4\pi t)^{3/2}} e^{-|x-y|^{2}/4t}\quad x,y\in \Pi,\ t>0.   \tag{3.5}
\end{align*}
\end{prop}

\vspace{5pt}

\begin{proof}
The pointwise upper bound (3.5) is known for an arbitrary domain. See \cite[Example 2.1.8]{Davies}. We prove the assertion (i). It suffices to show (3.4) for $f\in C^{\infty}_{c}$. Suppose that $f\geq 0$. Then, $u=e^{t\Delta_{D}}f$ is non-negative by a maximum principle. By multiplying $qu^{q-1}$ to the heat equation and integration by parts, we see that $\varphi=u^{q/2}$ satisfies 

\begin{align*}
\frac{d}{dt}\int_{\Pi}|\varphi|^{2}\dd x+\frac{4}{q'}\int_{\Pi}|\nabla \varphi|^{2}\dd x=0.
\end{align*}\\
Since the function $\varphi$ vanishes on $\partial\Pi$, we apply the Poincar\'e inequality in the cylinder $||\varphi||_{L^{2}}\leq C||\nabla \varphi||_{L^{2}} $\cite[6.30 THEOREM]{Ad} to estimate 

\begin{align*}
\frac{\dd}{\dd t} \int_{\Pi}|\varphi|^{2}\dd x\leq -\frac{4}{C^{2}q'}\int_{\Pi}|\varphi|^{2}\dd x.
\end{align*}\\
Thus the estimate (3.4) holds with the constant $M=4/C^{2}$. For general $f\in C^{\infty}_{c}$, we approximate the absolute value function as in the proof of Proposition 2.4 and obtain (3.4).
\end{proof}

\vspace{15pt}
Proposition 3.2 implies that the operator $-\Delta_{D}$ is invertible on $L^{q}$.

\vspace{15pt}

\begin{proof}[Proof of Lemma 3.1]
We prove (i). We set

\begin{align*}
\phi=\int_{0}^{\infty}e^{t \Delta_{D}}f\dd t\quad \textrm{for}\ f\in L^{q}.
\end{align*}\\
Since the heat semigroup is an analytic semigroup on $L^{q}$, it follows from (3.4) that 

\begin{align*}
||\phi||_{W^{1,q}}\leq C\Big(1+\frac{1}{\mu_q}\Big)||f||_{L^{q}}  \tag{3.6}
\end{align*}\\
with some constant $C$, independent of $q$. Since $-\Delta\phi=f$, by the elliptic regularity estimate \cite{ADN}, if follows that

\begin{align*}
||\nabla^{2}\phi||_{L^{q}}\leq C(||f||_{L^{q}}+||\phi||_{W^{1,q}})    \tag{3.7}
\end{align*}\\
We obtained (3.2). The uniqueness follows from a maximum principle.

We prove (ii). Since $\phi=(-\Delta_{D})^{-1}f=(-\Delta_{D})^{-1/2}(-\Delta_{D})^{-1/2}f$, we use a fractional power of the operator $-\Delta_{D}$. We set the domain $D(-\Delta_{D})$ by a space of all functions in $W^{2,q}$, vanishing on $\partial\Pi$. By estimates of pure imaginary powers of the operator \cite{Seeley71}, the domain of the fractional power $D((-\Delta_{D})^{1/2})$ is continuously embedded to  the Sobolev space $W^{1,q}$. Hence the operator $\partial (-\Delta_{D})^{-1/2}$ acts as a bounded operator on $L^{q}$.

It suffices to show that the fractional power $(-\Delta_{D})^{-1/2}$ acts as a  bounded operator from $L^{q}$ to $L^{p}$. We see that 

\begin{align*}
((-\Delta_{D})^{-1/2}f)(x)
=\int_{0}^{\infty}t^{-1/2}e^{t\Delta_{D}}f\dd t
=\int_{\Pi}f(y)\dd y\int_{0}^{\infty}t^{-1/2}K(x,y,t)\dd t.
\end{align*}\\
By (3.5), we have 

\begin{align*}
\int_{0}^{\infty}t^{-1/2}K(x,y,t)\dd t\leq \frac{C}{|x-y|^{2}}.
\end{align*}\\
Since the operator $f\longmapsto |x|^{-2}*f$ acts as a bounded operator from $L^{q}$ to $L^{p}$ for $1/p=1/q-1/3$ by the Hardy-Littlewood-Sobolev inequality \cite[p.354]{Stein93}, so is $(-\Delta_{D})^{-1/2}$. The proof is complete. 
\end{proof}

\vspace{15pt}

\subsection{Dependence of a constant}
The growth rate of the constant in (3.2) is at most linear as $q\to \infty$.

\vspace{15pt}

\begin{lem}
Let $q_0\in (3,\infty)$. There exists a constant $C$ such that 

\begin{align*}
||\phi||_{W^{2,q}(\Pi)}\leq C q||f||_{L^{q}\cap L^{q_0}(\Pi)}  \tag{3.8}
\end{align*}\\
holds for solutions of (3.1) for $f\in L^{q}\cap L^{q_0}(\Pi)$ and $q\in [q_{0},\infty)$, where 

\begin{align*}
||f||_{L^{q}\cap L^{q_0}}=\max\{||f||_{L^{q}}, ||f||_{L^{q_0}}\}.
\end{align*}

\end{lem}

\vspace{15pt}
We prove Lemma 3.3 by a cut-off function argument.

\vspace{15pt}

\begin{prop}
Let $G$ be a smoothly bounded domain in $\mathbb{R}^{3}$. Let $q_0\in (3,\infty)$. There exists a constant $C$ such that 

\begin{align*}
||\phi||_{W^{2,q}(G)}\leq C q||f||_{L^{q}(G)}  \tag{3.9}
\end{align*}\\
holds for solutions of (3.1) for $f\in L^{q}(G)$ and $q\in [q_0,\infty)$.
\end{prop}

\vspace{5pt}

\begin{proof}
The assertion is proved in \cite[Corollary 1]{Yudovich62} for general elliptic operators and $n$-dimensional bounded domains.
\end{proof}

\vspace{5pt}

\begin{proof}[Proof of Lemma 3.3]
Let $\{\varphi_{j}\}_{j=-\infty}^{\infty}\subset C^{\infty}_{c}(\mathbb{R})$ be a partition of the unity such that $0\leq \varphi_j\leq 1$, $\textrm{spt}\ \varphi_j\subset [j-1,j+1]$ and $\sum_{j=-\infty}^{\infty}\varphi_j(x_3)=1$, $x_3\in \mathbb{R}$. Let $\phi\in W^{2,q}(\Pi)$ be a solution of (3.1) for $f\in L^{q}\cap L^{q_0}(\Pi)$. We set $\phi_j=\phi\varphi_j$ and observe that 

\begin{align*}
-\Delta \phi_j&=f_j\quad \textrm{in}\ G_j,\\
\phi_j&=0\quad \textrm{on}\ \partial G_j,
\end{align*}\\
for $G_j=D\times (j-1,j+1)$ and $f_j=f\varphi_j-2\nabla \phi\cdot \nabla \varphi_j-\phi\Delta \varphi_j$. We take a smooth bounded domain $\tilde{G}_{j}$ such that $G_j\subset \tilde{G}_{j}\subset D\times [j-2,j+2]$ and apply (3.9) to estimate 

\begin{align*}
||\phi_j||_{W^{2,q}(\tilde{G}_j)}\leq Cq ||f_{j}||_{L^{q}(\tilde{G}_{j})}
\end{align*}\\
for $q\in [q_0,\infty)$ with some constant $C$, independent of $j$ and $q$. It follows that

\begin{align*}
||\nabla^{2}\phi \varphi_{j}||_{L^{q}(\Pi)}\leq Cq(||f||_{L^{q}(G_j)}+||\phi||_{W^{1,q}(G_j)} ).
\end{align*}\\
By summing over $j$, we obtain 

\begin{align*}
||\nabla^{2}\phi ||_{L^{q}(\Pi)}\leq Cq(||f||_{L^{q}(\Pi)}+||\phi||_{W^{1,q}(\Pi)} ).   \tag{3.10}
\end{align*}\\
We estimate the lower order term of $\phi$. By Lemma 3.1(i), we have $||\phi||_{W^{2,q_0}}\leq C||f||_{L^{q_0}}$. In particular, $||\phi||_{W^{1,\infty}}\leq C||f||_{L^{q_0}}$ by the Sobolev inequality. Applying the H\"older inequality implies that 

\begin{align*}
||\phi||_{W^{1,q}(\Pi)}\leq C||f||_{L^{q_0}(\Pi)}  \tag{3.11}
\end{align*}\\
for $q\in [q_0,\infty)$ with some constant $C$, independent of $q$. The estimate (3.8) follows from (3.10) and (3.11).
\end{proof}

\vspace{15pt}

\subsection{Stream functions}
We shall give a Biot-Savart law for axisymmetric vector fields without swirl. We see that a smooth axisymmetric solenoidal vector field without swirl $u=u^{r}e_{r}+u^{z}e_{z}$ in $\Pi$ satisfies  

\begin{align*}
\partial_z(ru^{z})+\partial_r(ru^{r})&=0\quad (z,r)\in \mathbb{R}\times (0,1),\\
ru^{r}&=0\quad \textrm{on}\ \{r=0,1\}.
\end{align*}\\
Since $(ru^{z},ru^{r})$ is regarded as a solenoidal vector field in the two-dimensional layer $\mathbb{R}\times (0,1)$, there exists a stream function $\psi(r,z)$ such that 

\begin{align*}
ru^{z}=\frac{\partial \psi}{\partial r},\quad ru^{r}=-\frac{\partial \psi}{\partial z}. 
\end{align*}\\
Since $\psi$ is constant on the boundary, we may assume that $\psi=0$ on $\{r=1\}$. Since $\phi=(\psi/r)e_{\theta}$ satisfies 

\begin{align*}
\D\ \phi=0,\ \nabla \times \phi=u\quad \textrm{in}\ \Pi,\quad \phi=0\quad \textrm{on}\ \partial\Pi,
\end{align*}\\
we see that $-\Delta \phi=\nabla \times u$. Since the Laplace operator $-\Delta_{D}$ is invertible, the stream function is represented by $\phi=(-\Delta_{D})^{-1}\nabla \times u$.

\vspace{15pt}

\begin{lem}
(i) Let $u$ be an axisymmetric vector field without swirl in $L^{p}_{\sigma}$ such that $\nabla \times u\in L^{q}$ for $q\in (1, 3)$ and $1/p=1/q-1/3$. Then, 

\begin{align*}
u=\nabla \times (-\Delta_{D})^{-1} (\nabla\times u).    \tag{3.12}
\end{align*}\\

\noindent 
(ii) The estimates

\begin{align*}
||u||_{L^{p}}+||\nabla u||_{L^{q}}&\leq C_1||\nabla \times u||_{L^{q}},  \tag{3.13}\\
||\nabla u||_{L^{r}}&\leq C_2||\nabla \times u||_{L^{r}},\quad 1<r<\infty, \tag{3.14} \\
||\nabla u||_{L^{r}}&\leq C_3r||\nabla \times u||_{L^{r}\cap L^{r_0}},\quad 3<r_0< r<\infty, \tag{3.15} 
\end{align*}\\
hold with some constants $C_1-C_3$. The constant $C_3$ is independent of $r$.
\end{lem}

\vspace{15pt}
It suffices to show (3.12). The assertion (ii) follows from Lemmas 3.1 and 3.3.
\vspace{15pt}

\begin{prop}
Let $w$ be an axisymmetric vector field without swirl in $L^{p}$ for $p\in (1,\infty)$. Assume that 

\begin{align*}
\D\ w=0,\ \nabla \times w=0\quad \textrm{in}\ \Pi,\quad w\cdot n=0\quad \textrm{on}\ \partial\Pi.   
\end{align*}\\
Then, $w\equiv 0$.
\end{prop}

\vspace{5pt}

\begin{proof}
Since $w=w^{r}e_{r}+w^{z}e_{z}$ is a harmonic vector field in $\Pi$ and $\Delta=\partial_r^{2}+r^{-1}\partial_r+r^{-2}\partial_{\theta}^{2}+\partial_{z}^{2}$ by the cylindrical coordinate, we see that 

\begin{align*}
0=\Delta w&=\Delta (w^{r}e_{r})+\Delta (w^{z}e_{z}) \\
&=\left\{\left(\partial_r^{2}+r^{-1}\partial_r-r^{-2}+\partial_{z}^{2}\right)w^{r}\right\}e_{r}+(\Delta w^{z})e_{z}\\
&=\left\{\left(\Delta-r^{-2}\right)w^{r}\right\}e_{r}+(\Delta w^{z})e_{z}.
\end{align*}\\
Hence, $(\Delta-r^{-2})w^{r}=0$ and $\Delta w^{z}=0$. By $\Delta (w^{r}e_r)=\{(\Delta-r^{-2})w^{r}\}e_r=0$, $w^{r}e_{r}$ is harmonic in $\Pi$. Since $w^{r}$ vanishes on the boundary and the operator $-\Delta_{D}$ is invertible on $L^{p}$, we see that $w^{r}\equiv 0$. By the divergence-free condition $\partial_r w^{r}+w^{r}/r+\partial_z w^{z}=0$ and a decay condition $w^{z}\in L^{p}$, we have $w^{z}\equiv 0$.
\end{proof}

\vspace{15pt}

\begin{proof}[Proof of Lemma 3.5]
We set 

\begin{align*}
\tilde{\phi}=(-\Delta_{D})^{-1}(\nabla \times u),\quad \tilde{u}=\nabla \times \tilde{\phi}.
\end{align*}\\
Since $u$ is axisymmetric without swirl, $\tilde{\phi}$ is axisymmetric and $\tilde{\phi}=\tilde{\phi}^{\theta} e_{\theta}$. Since $\tilde{\phi}$ satisfies 

\begin{align*}
\D\ \tilde{\phi}=0,\quad -\Delta \tilde{\phi}=\nabla \times u\quad \textrm{in}\ \Pi,\qquad \tilde{\phi}=0\quad \textrm{on}\ \partial\Pi,
\end{align*}\\
it follows that 

\begin{align*}
\D\ \tilde{u}=0,\quad \nabla \times \tilde{u}=\nabla \times u\quad \textrm{in}\ \Pi,\quad \tilde{u}\cdot n=0\quad \textrm{on}\ \partial\Pi.
\end{align*}\\
Applying Proposition 3.6 for $w=u-\tilde{u}$ implies $u\equiv \tilde{u}$. We proved (3.12). 
\end{proof}

\vspace{15pt}

\section{Vanishing viscosity limits}

\vspace{15pt}
We prove Theorem 1.2 (i) and (ii). When the initial vorticity satisfies $\omega^{\theta}_{0}/r\in L^{q}$ for $q\in [3/2,3)$, the initial velocity belongs to $L^{p}$ for $p\in [3,\infty)$ by the Biot-Savart law and a global-in-time unique solution $u=u_{\nu}$ of (1.1) exists by Theorem 1.1. We use the vorticity estimate (2.3) and construct global weak solutions of the Euler equations by sending $\nu\to0$. In the subsequent section, we prove H\"older continuity of weak solutions.

\vspace{15pt}

\subsection{Convergence to a limit}

We first derive a priori estimates independent of the viscosity $\nu>0$.

\vspace{15pt}

\begin{lem}
(i) Let $u_0\in L^{p}_{\sigma}$ be an axisymmetric vector field without swirl such that $\omega^{\theta}_{0}/r\in L^{q}$ for $q\in [3/2,3)$ and $1/p=1/q-1/3$. Let $u_{\nu}\in C([0,\infty); L^{p})\cap C^{\infty}(\overline{\Pi}\times (0,\infty))$ be a solution of (1.1) for $u_0$ in Theorem 1.1. There exits a constant $C$ such that 

\begin{align*}
||u_\nu||_{L^{p}(\Pi)}+||\nabla u_\nu||_{L^{q}(\Pi)}\leq C\Big\|\frac{\omega^{\theta}_{0}}{r}\Big\|_{L^{q}(\Pi)} \quad t\geq 0,\ \nu>0.  \tag{4.1}
\end{align*}\\
Moreover, for each bounded domain $G\subset \Pi$, there exists a constant $C'$ such that 

\begin{align*}
||\partial_t u_{\nu}||_{W^{-1,q}(G)}\leq C'\Big\|\frac{\omega^{\theta}_{0}}{r}\Big\|_{L^{q}(\Pi)}\Big(\nu+\Big\|\frac{\omega^{\theta}_{0}}{r}\Big\|_{L^{q}(\Pi)}\Big) \quad t\geq 0,\ \nu>0, \tag{4.2}
\end{align*}\\
where $W^{-1,q}$ denotes the dual space of $W^{1,q'}_{0}$ and $q'$ is the conjugate exponent to $q$.

\noindent 
(ii) If $\omega^{\theta}_{0}/r\in L^{s}$ for $s\in (3,\infty)$, then

\begin{align*}
||\nabla u_\nu||_{L^{s}(\Pi)}\leq C\Big\|\frac{\omega^{\theta}_{0}}{r}\Big\|_{L^{s}(\Pi)} \quad t\geq 0,\ \nu>0.  \tag{4.3}
\end{align*}\\

\noindent 
(iii) If $\omega^{\theta}_{0}/r\in L^{\infty}$, then

\begin{align*}
||\nabla\times  u_\nu||_{L^{\infty}(\Pi)}\leq \Big\|\frac{\omega^{\theta}_{0}}{r}\Big\|_{L^{\infty}(\Pi)} \quad t\geq 0,\ \nu>0.  \tag{4.4}
\end{align*}
\end{lem}

\vspace{15pt}

\begin{proof}
Since $\omega^{\theta}_{0}/r\in L^{q}$, applying Lemma 2.3 implies the vorticity estimate

\begin{align*}
\Big\|\frac{\omega^{\theta}_{\nu}}{r}\Big\|_{L^{q}}\leq \Big\|\frac{\omega^{\theta}_{0}}{r}\Big\|_{L^{q}}.
\end{align*}\\
Since $r<1$, the $L^{q}$-norm of the vorticity $\nabla \times u_{\nu}$ is bounded. By the estimate of the Biot-Savart law (3.13), we obtain (4.1). The estimates (4.3) and (4.4) follow in the same way.

We prove (4.2). We take an arbitrary $\varphi\in C^{\infty}_{c}(G)$ and consider its zero extension to $\Pi\backslash \overline{G}$ (denoted by $\varphi$). We set $f=\mathbb{P}\varphi$ by the Helmholtz projection operator $\mathbb{P}$. By a higher regularity estimate of the Helmholtz projection operator \cite[Theorem 6]{ST98}, we see that $f\in C^{\infty}(\overline{\Pi})$ and 

\begin{align*}
|| f||_{W^{1,s}}= || \mathbb{P}\varphi||_{W^{1,s}}\leq C_s|| \varphi||_{W^{1,s}}
\end{align*}\\
with some constant $C_s$ and $s\in(1,\infty)$. By multiplying $f$ by (1.1) and integration by parts, we see that 

\begin{align*}
\int_{\Pi}\partial_t u_{\nu}\cdot f\dd x 
&=\nu\int_{\Pi} \Delta u_{\nu}\cdot f\dd x
-\int_{\Pi}(u_{\nu}\cdot \nabla u_{\nu})\cdot f\dd x\\
&=-\nu \int_{\Pi}\nabla \times u_{\nu}\cdot \nabla \times f\dd x
+\int_{\Pi}u_{\nu}u_{\nu}: \nabla f\dd x.
\end{align*}\\
By $\D\ u=0$, the left-hand side equals to the integral of $\partial_{t} u\cdot \varphi$ in $G$. By applying the estimate of the Helmholtz projection, we obtain 

\begin{align*}
\Bigg|\int_{G}\partial_t u_{\nu}\cdot \varphi\dd x\Bigg|
\leq C\Bigg(\nu ||\nabla u_{\nu}||_{L^{q}(\Pi)}|| \varphi||_{W^{1,q'}(G)}+||u_{\nu}||_{L^{p}(\Pi)}^{2}|| \varphi||_{W^{1,p/(p-2)}(G)}   \Bigg),
\end{align*}\\
Since $p/(p-2)\leq q'$, the norms of $\varphi$ are estimated by the $W^{1,q'}$-norm of $\varphi$ in $G$. By (4.1), we obtain (4.2).
\end{proof}

\vspace{15pt}

We apply the estimates (4.1) and (4.2) in order to extract a subsequence of $\{u_{\nu}\}$. We recall an abstract compactness theorem in \cite[Chapter III, Theorem 2.1]{Te}.

\vspace{15pt}

\begin{prop}
(i) Let $X_0$, $X$ and $X_1$ be Banach spaces such that $X_0\subset X\subset X_1$ with continuous injections, $X_0$ and $X_1$ are reflexive and the injection $X_0\subset X$ is compact. For $T\in (0,\infty)$ and $s\in (1,\infty)$, set the Banach space

\begin{align*}
Y=\{u\in L^{s}(0,T; X_0)\ |\ \partial_t u\in L^{s}(0,T; X_1)  \},
\end{align*}\\
equipped with the norm $||u||_{Y}=||u||_{L^{s}(0,T; X_0)}+||\partial_t u||_{L^{s}(0,T; X_1)}$. Then, the injection $Y\subset L^{s}(0,T; X)$ is compact.
\end{prop}

\vspace{15pt}

\begin{proof}[Proof of Theorem 1.2 (i)]
For an arbitrary bounded domain $G\subset \Pi$, we set $X_0=W^{1,q}(G)$, $X=L^{q}(G)$ and $X_1=W^{-1,q}(G)$. Since $\{u_\nu\}$ is a bounded sequence in $Y$ by (4.1) and (4.2), we apply Proposition 4.2 to get a subsequence (still denoted by $u_\nu$) that converges to a limit $u$ in $L^{s}(0,T; L^{q}(G))$. By choosing a subsequence, we may assume that $u_\nu$ converges to $u$ in $L^{s}(0,T; L^{q}(G))$ for arbitrary $G\subset \Pi$ and $T>0$ and 

\begin{align*}
u_\nu \to u\quad \textrm{a.e. in}\ \Pi\times (0,\infty).
\end{align*}\\
We take an arbitrary $\varphi\in C^{1}_{c}(\overline{\Pi}\times [0,\infty))$ such that $\D\ \varphi=0$ in $\Pi$ and $\varphi\cdot n=0$ on $\partial\Pi$. By multiplying $\varphi$ by (1.1) and integration by parts, we see that 

\begin{align*}
\int_{0}^{\infty}\int_{\Pi} (u_{\nu}\cdot \partial_t \varphi -\nu\nabla u_{\nu}\cdot \nabla \varphi +u_{\nu}u_{\nu}:\nabla \varphi)\dd x\dd t
=-\int_{\Pi}u_0\cdot \varphi_0\dd x   \tag{4.5}
\end{align*}\\
Note that the integral of $\partial_n u_{\nu}\cdot \varphi$ on $\partial\Pi$ vanishes since $\partial_r u^{z}_{\nu}=0$ and $\varphi\cdot n=0$ on $\partial\Pi$. The first term converges to the integral of $u\cdot \partial_t \varphi$ and the second term vanishes by (4.1). We take a bounded domain $G\subset \Pi$ and $T>0$ such that $\textrm{spt}\ \varphi\subset \overline{G}\times [0,T]$. Since $u_{\nu}$ converges to $u$ a.e. in $\Pi\times (0,\infty)$, by Egoroff's theorem (e.g., \cite[1.2]{EG}) for an arbitrary $\varepsilon>0$, there exists a measurable set $E\subset G\times(0,T)$ such that $|G\times (0,T)\backslash E|\leq \varepsilon$ and 

\begin{align*}
u_{\nu}\to u\quad \textrm{uniformly on}\ E.
\end{align*}\\
Here, $|\cdot |$ denotes the Lebesgue measure. We set $F=G\times (0,T)\backslash E$. It follows from (4.1) that 

\begin{align*}
\Bigg|\int\hspace{-5pt}\int_{F} u_{\nu}u_{\nu}: \nabla \varphi\dd x\dd t\Bigg|
\leq ||u_\nu||_{L^{p}(F)}^{2}||\varphi||_{L^{p/(p-2)}(F)}
\leq C\varepsilon^{1-2/p},
\end{align*}\\
with some constant $C$, independent of $\nu$ and $\varepsilon$. Since $\varepsilon>0$ is arbitrary and $u_\nu$ converges to $u$ uniformly on $E$, the third term of (4.5) converges to the integral of $uu:\nabla \varphi$. Thus sending $\nu\to 0$ to (4.5) yields (1.5). Since the estimate (4.1) is inherited to the limit $u$, we see that $u\in L^{\infty}(0,\infty; L^{p})$ and $\nabla u\in L^{\infty}(0,\infty; L^{q})$.

We show the weak continuity $u\in BC_{w}([0,\infty); L^{p})$. We take an arbitrary $\varphi\in C_{c}^{\infty}(\Pi)$ and $\eta\in C^{1}[0,\infty)$. By multiplying $\varphi\eta$ by (1.1) and integration by parts as we did in the proof of Lemma 4.1, we obtain the estimate 

\begin{align*}
&\Bigg|\int_{0}^{\infty}\Bigg(\int_{\Pi}u_{\nu}(x,t)\cdot \varphi(x)\dd x\Bigg)\dot{\eta}(t)\dd t\Bigg|\\
&\leq C\Big\|\frac{\omega^{\theta}_{0}}{r}\Big\|_{L^{q}(\Pi)}\Big(\nu || \varphi||_{W^{1,q'}(\Pi)}+ \Big\|\frac{\omega^{\theta}_{0}}{r}\Big\|_{L^{q}(\Pi)} || \varphi||_{W^{1,p/(p-2)}(\Pi)}\Big)\int_{0}^{\infty}\eta(s)\dd s
\end{align*}\\
Since the left-hand side converges, the integral of $u\varphi$ in $\Pi$ is weakly differentiable as a function of time. By sending $\nu \to 0$ and the duality, we obtain 

\begin{align*}
\Bigg|\frac{\dd }{\dd t}\int_{\Pi}u(x,t)\cdot \varphi(x)\dd x\Bigg|
\leq C\Big\|\frac{\omega^{\theta}_{
0}}{r}\Big\|_{L^{q}(\Pi)}|| \varphi||_{W^{1,p/(p-2)}(\Pi)}\quad \textrm{a.e.}\ t>0.
\end{align*}\\
Hence for $s\in [0,\infty)$, we have

\begin{align*}
\int_{\Pi}u(x,t)\cdot \varphi(x)\dd x\to \int_{\Pi}u(x,s)\cdot \varphi(x)\dd x \quad \textrm{as}\ t\to s.
\end{align*}\\
By $u\in L^{\infty}(0,\infty; L^{p})$ and the density, the above convergence holds for all $\varphi\in L^{p'}$. Thus $u\in BC_{w}([0,\infty); L^{p})$. The weak continuity of $\nabla u$ on $L^{q}$ follows from that of $u$ on $L^{p}$. We proved the assertion (i).
\end{proof}

\vspace{15pt}

\subsection{Regularity of weak solutions}

We prove Theorem 1.2 (ii). We use the Poincar\'e inequality.

\vspace{15pt}

\begin{prop}
Let $u\in C(\overline{\Pi})$ be an axisymmetric vector field without swirl such that $\D\ u=0$ in $\Pi$, $u\cdot n=0$ on $\partial\Pi$ and $u(x)\to0$ as $|x|\to\infty$. Assume that $\nabla u\in L^{s}(\Pi)$ for $s\in (1,\infty)$. Then, the estimate 

\begin{align*}
||u||_{L^{s}(\Pi)}\leq C||\nabla u||_{L^{s}(\Pi)}  \tag{4.6}
\end{align*}\\
holds with some constant $C$.
\end{prop}

\vspace{5pt}

\begin{proof}
Since the radial component $u^{r}$ vanishes on $\partial\Pi$ by $u\cdot n=0$, we apply the Poincar\'e inequality \cite{Ad} to estimate 

\begin{align*}
||u^{r}||_{L^{s}(\Pi)}\leq C ||\nabla u^{r}||_{L^{s}(\Pi)}. 
\end{align*}\\
We estimate $u^{z}$. For arbitrary $z_1, z_2\in \mathbb{R}$, we set $G=D\times (z_1,z_2)$. Since $\D\ u=0$, it follows that 

\begin{align*}
0=\int_{G}\D\ u\dd x=\int_{D}u^{z}(r,z_2)\dd {\mathcal{H}}
-\int_{D}u^{z}(r,z_1)\dd {\mathcal{H}}.
\end{align*}\\
Since $u$ decays as $|z_2|\to\infty$, we see that the flux on $D$ is zero, i.e., 

\begin{align*}
\int_{D}u^{z}(r,z_1)\dd {\mathcal{H}}=0.
\end{align*}\\
We apply the Poincar\'e inequality \cite{E} to estimate

\begin{align*}
||u^{z}||_{L^{s}(D)}(z_1)\leq C||\nabla_{h} u^z||_{L^{s}(D)}(z_1),
\end{align*}\\
where $\nabla_h$ denotes the gradient for the horizontal variable $x_h=(x_1,x_2)$. By integrating for $z_1\in \mathbb{R}$, we obtain (4.6).
\end{proof}

\vspace{15pt}

\begin{proof}[Proof of Theorem 1.2 (ii)]
If $\omega^{\theta}_{0}/r\in L^{s}$ for $s\in (3,\infty)$, the limit $u\in BC_{w}([0,\infty); L^{p})$ satisfies $\nabla u\in L^{\infty}(0,\infty; L^{s})$ by (4.3). Thus $u(\cdot ,t)$ is H\"older continuous in $\overline{\Pi}$ and decaying as $|x|\to\infty$. We apply the Poincar\'e inequality (4.6) to see that $u\in L^{\infty}(0,\infty; W^{1,s})$ and $u\cdot \nabla u\in L^{\infty}(0,\infty; L^{s})$ by the Sobolev inequality. By integration by parts, it follows from (1.3) that 

\begin{align*}
\int_{0}^{\infty}\int_{\Pi}u\cdot \varphi\dot{\eta}\dd x\dd t
=\int_{0}^{\infty}\int_{\Pi} (u\cdot \nabla u)\cdot \varphi\eta\dd x\dd t \tag{4.7}
\end{align*}\\
for all $\varphi\in L^{s'}_{\sigma}$ and $\eta\in C^{\infty}_{c}(0,\infty)$, where $s'$ is the conjugate exponent to $s$. By the boundedness of the Helmholtz projection on $L^{s'}$ and a duality, we see that $\partial_t u\in L^{\infty}(0,\infty; L^{s})$ and $u\in  BC([0,\infty); L^{s})$. The equation (1.4) follows from (4.7) by integration by parts. Since $\nabla u$ is bounded and $u$ is continuous on $L^{s}$, $\nabla u$ is weakly continuous on $L^{s}$. We proved the assertion (ii).

If in addition that $\omega^{\theta}_{0}/r\in L^{\infty}$, the limit satisfies $\nabla \times u\in BC_{w}([0,\infty); L^{\infty})$ by (4.4).
\end{proof}

\vspace{15pt}

\begin{rems}
(i) The equation (1.4) is written as

\begin{align*}
\partial_t u+u\cdot \nabla u+\nabla p=0\qquad \textrm{on}\ L^{s}\quad \textrm{for a.e.}\ t>0,
\end{align*}\\
by the associated pressure $\nabla p=-(I-\mathbb{P})u\cdot \nabla u\in L^{\infty}(0,\infty; L^{s})$. 

\noindent
(ii) The weak solutions in Theorem 1.2 are with finite energy for $q\in [3/2,2]$. In fact, by the global estimate (4.1) and applying the Poincar\'e inequality (4.6) for global-in-time solutions $u=u_{\nu}$ of (1.1) for $\omega^{\theta}_{0}/r\in L^{q}$ in Theorem 1.1, we see that 

\begin{align*}
||u_{\nu}||_{L^{p}}+||u_{\nu}||_{L^{q}}\leq C||\nabla u_{\nu}||_{L^{q}}\leq C'\Big\|\frac{\omega^{\theta}_{0}}{r}\Big\|_{L^{q}}\qquad t\geq 0,\ \nu>0,
\end{align*}\\
for $p\in [3,6]$ satisfying $1/p=1/q-1/3$. By the H\"older inequality, the solutions are uniformly bounded in $L^{r}$ for all $r\in [2,3]$ and the limit belongs to the same space. 
\end{rems}

\vspace{15pt}

\section{Uniqueness}

\vspace{15pt}

We prove Theorem 1.2 (iii). It remains to show the uniqueness. Since the weak solutions are with infinite energy for $q\in (2,3)$, we estimate a local energy of two weak solutions in the cylinder by using a cut-off function $\theta_{R}$. We then send $R\to\infty$ and prove the uniqueness by using the growth bound of the $L^{r}$-norm (3.15). To this end, we show decay properties of weak solutions.

\vspace{15pt}

\subsection{Decay properties of weak solutions}

We use the Poincar\'e inequality (4.6) and deduce decay properties of velocity as $|x_3|\to\infty$.

\vspace{15pt}

\begin{prop}
The weak solutions $(u,p)$ in Theorem 1.2 (iii) satisfy

\begin{align*}
u,\nabla u,\partial_t u,\nabla p\in L^{\infty}(0,\infty; L^{q}).  \tag{5.1}
\end{align*}
\end{prop}

\vspace{5pt}

\begin{proof}
Since $u(\cdot,t)$ is bounded and H\"older continuous in $\overline{\Pi}$ and $\nabla u\in BC_{w}([0,\infty); L^{q})$, we see that $u\cdot \nabla u\in L^{\infty}(0,\infty; L^{q})$. Thus, $\partial_t u$ and $\nabla p=-(I-\mathbb{P})u\cdot \nabla u$ belong to $L^{\infty}(0,\infty; L^{q})$ by (1.4). By the Poincar\'e inequality (4.6), $u\in L^{\infty}(0,\infty; L^{q})$ follows. 
\end{proof}

\vspace{15pt}

We estimate the pressure $p$ as $|x_3|\to\infty$. We set

\begin{equation*}
\begin{aligned}
&\tilde{p}(x_h,x_3,t)=p(x_h,x_3, t)-\hat{p}(x_3, t),\\
&\hat{p}(x_3,t)=\frac{1}{|D|}\int_{D}p(x_h,x_3,t)\dd x_h.
\end{aligned}
\tag{5.2}
\end{equation*}

\vspace{15pt}

\begin{prop}
\begin{align*}
&\tilde{p}\in L^{\infty}(0,\infty; L^{q}),   \tag{5.3} \\
&|\hat{p}(x_3,t)|\leq C(1+|x_3|)^{1/3}\quad x_3\in \mathbb{R},\ t>0,  \tag{5.4}
\end{align*}\\
with some constant $C$.
\end{prop}

\vspace{5pt}

\begin{proof}
The property (5.3) follows from (5.1) by applying the Poincar\'e inequality on $D$. We show (5.4). We integrate the vertical component of (1.2) on $D$ to see that 

\begin{align*}
\frac{\partial}{\partial t} \int_{D}u^{z}\dd x_h
+\int_{D}u\cdot \nabla u^{z}\dd x_h+\frac{\partial}{\partial z}\int_{D}p\dd x_h=0.
\end{align*}\\
Since the flux of $u$ on $D$ is zero as we have seen in the proof of Proposition 4.3, the first term vanishes. We integrate the equation by the vertical variable between $(0,z)$ to get

\begin{align*}
\int_{D}p(r,z,t)\dd x_h
=\int_{D}p(r,0,t)\dd x_h
-\int_{0}^{z}\int_{D}u\cdot \nabla u^{z}\dd x.
\end{align*}\\
We observe that $u\in L^{\infty}(0,\infty; W^{1,r})$ for $r\in [q,3)$ by (5.1) and Theorem 1.2 (ii). Since $u\cdot \nabla u\in L^{\infty}(0,\infty; L^{r/2})$ for $r\geq 2$, we apply the H\"older inequality to estimate 

\begin{align*}
\Bigg|\int_{0}^{z}\int_{D}u\cdot \nabla u^{z}\dd x\Bigg|
\leq |D\times (0,z)|^{1-2/r}||u\cdot \nabla u||_{L^{r/2}}
\leq C|z|^{1-2/r}\qquad t>0.
\end{align*}\\
Since $r\in [2,3)$ and $1-2/r<1/3$, we obtain (5.4).
\end{proof}

\vspace{15pt}
We use the growth bound for the $L^{r}$-norm of $\nabla u$ as $r\to\infty$.

\vspace{15pt}

\begin{prop}
\begin{align*}
&\nabla u\in BC_{w}([0,\infty); L^{r})\quad r\in (3,\infty),   \\
&\nabla\times  u\in BC_{w}([0,\infty); L^{\infty}),   \\
&||\nabla u||_{L^{r}}\leq Cr\quad r>3,\ t\geq 0,  
\end{align*}\\
with some constant $C$, independent of $r$.
\end{prop}

\vspace{5pt}

\begin{proof}
By construction, the weak solution $u$ satisfies 

\begin{align*}
||\nabla\times u||_{L^{r}}\leq \Big\|\frac{\omega^{\theta}_{0}}{r}\Big\|_{L^{r}}\qquad t\geq 0.
\end{align*}\\
for all $r\in (3,\infty)$. We fix $r\in (3,\infty)$ and take $r_0\in (3,r)$. It follows from (3.15) that 

\begin{align*}
||\nabla u||_{L^{r}}\leq C r ||\nabla \times u||_{L^{r}\cap L^{r_0}}
\leq C r \Big\|\frac{\omega^{\theta}_{0}}{r}\Big\|_{L^{r}\cap L^{r_0}}.
\end{align*}\\
Since the $L^{r}$-norm of $\omega^{\theta}_{0}/r$ is uniformly bounded for all $r>3$ by $\omega^{\theta}_{0}/r\in L^{\infty}$, we obtain the desired estimate for $\nabla u$.
\end{proof}

\vspace{15pt}

\subsection{Local energy estimates}

\vspace{15pt}

We now prove the uniqueness. Let $(u_1,p_1)$ and $(u_2,p_2)$ be two weak solutions to (1.2) in Theorem 1.2 (iii) for the same initial data. Then,  $w=u_1-u_2$ and $\pi=p_1-p_2$ satisfy

\begin{equation*}
\begin{aligned}
\partial_t w+u_1\cdot \nabla w+w\cdot \nabla u_2+\nabla \pi=0,\quad \D\ w&=0\qquad \textrm{in}\ \Pi\times (0,\infty),\\
w\cdot n&=0\qquad \textrm{on}\ \partial\Pi\times (0,\infty), \\
w&=0\qquad \textrm{on}\ \Pi\times\{t=0\}.
\end{aligned}
\tag{5.5}
\end{equation*}\\
Let $\theta\in C^{\infty}_{c}[0,\infty)$ be a smooth monotone non-increasing function such that $\theta\equiv 1$ in $[0,1]$ and $\theta\equiv 0$ in $[2,\infty)$. We set $\theta_R(x_3)=\theta(|x_3|/R)$ for $R\geq 1$ so that $\theta_{R}\equiv 1$ in $[0,R]$, $\theta_{R}\equiv 0$ in $[2R,\infty)$, $||\partial_{x_3}\theta_{R}||_{\infty}\leq C/R$ and $\textrm{spt}\ \partial_{x_3}\theta_{R}\subset I_R$ for $I_{R}=[R,2R]$. By multiplying $2w\theta_{R}$ by $(1.2)$ and integration by parts, we see that 

\begin{align*}
\frac{\dd}{\dd t}\int_{\Pi}|w|^{2}\theta_{R}\dd x
+2\int_{\Pi}(w\cdot \nabla u_2)\cdot w\theta_{R}\dd x
-\int_{\Pi}u_1|w|^{2}\cdot \nabla \theta_{R}\dd x
-2\int_{\Pi}\pi w\cdot \nabla \theta_{R}\dd x=0.   \tag{5.6}
\end{align*}\\
We set 

\begin{align*}
\phi_{R}(t)=\int_{\Pi}|w|^{2}(x,t)\theta_{R}(x_3)\dd x.
\end{align*}\\
By Theorem 1.2 (ii), the function $\phi_R\in C[0,\infty)$ is differentiable for a.e. $t>0$ and satisfies $\phi_R(0)=0$. We estimate errors in the cut-off procedure. 

\vspace{15pt}

\begin{prop}
There exists a constant $C=C(R)$ such that 

\begin{align*}
\Bigg|\int_{\Pi}u_1|w|^{2}\cdot \nabla \theta_{R}\dd x\Bigg|+
\Bigg|2\int_{\Pi}\pi w\cdot \nabla \theta_{R}\dd x\Bigg|\leq C\quad t>0.  \tag{5.7}
\end{align*}\\
The constant $C(R)$ converges to zero as $R\to\infty$ for each $t>0$.
\end{prop}

\vspace{5pt}

\begin{proof}
Since $u_1\in L^{\infty}(0,\infty; L^{q})$ and $w\in L^{\infty}(\Pi\times (0,\infty))$ by Proposition 5.1 and Theorem 1.2 (ii), applying the H\"older inequality yields

\begin{align*}
\Bigg|\int_{\Pi}u_1|w|^{2}\cdot \nabla \theta_{R}\dd x\Bigg|
\leq \frac{C}{R}\int_{D\times I_R}|u_1|\dd x 
\leq \frac{C'}{R}|D\times I_R|^{1/q'}||u_1||_{L^{q}}
\leq \frac{C''}{R^{1/q}}.
\end{align*}\\
We next estimate the second term of (5.7). We set $\pi=\tilde{\pi}+\hat{\pi}$ by (5.2). Since $\tilde{\pi}\in L^{\infty}(0,\infty; L^{q})$ by (5.3), it follows that 

\begin{align*}
\Bigg|\int_{\Pi}\tilde{\pi}w\cdot \nabla \theta_R\dd x\Bigg|
\leq \frac{C}{R}\int_{D\times I_R}|\tilde{\pi}|\dd x\leq \frac{C'}{R^{1/q}}.
\end{align*}\\
It follows from (5.4) that 

\begin{align*}
\Bigg|\int_{\Pi}\hat{\pi}w\cdot \nabla \theta_R\dd x\Bigg|
\leq \frac{C}{R^{2/3}}\int_{D\times I_R}|w|\dd x
\leq \frac{C}{R^{2/3}}|D\times I_R|^{1/q'}||w||_{L^{q}}
\leq \frac{C'}{R^{2/3-1/q'}}.
\end{align*}\\
Since $2/3-1/q'>0$ for $q\in [3/2,3)$, the right-hand side converges to zero as $R\to\infty$. 
\end{proof}

\vspace{15pt}

\begin{proof}[Proof of Theorem 1.2 (iii)]
By Proposition 5.3, there exist constants $M_1$ and $M_2$ such that 

\begin{align*}
||w||_{L^{\infty}}&\leq M_1,\\
||\nabla u_2||_{L^{r}}&\leq M_2 r\qquad r>3,\ t\geq 0.
\end{align*}\\
For an arbitrary $\delta\in (0,2/3)$, we set $r=2/\delta$. We apply the H\"older inequality with the conjugate exponent $r'=2/(2-\delta)$ to see that 

\begin{align*}
\Bigg|2\int_{\Pi}(w\cdot \nabla u_2)\cdot w\theta_{R}\dd x\Bigg|
&\leq 2\int_{\Pi}|\nabla u_2|(|w|\theta_{R}^{1/2})^{2}\dd x\\
&\leq 2M_1^{\delta}\int_{\Pi}|\nabla u_2|(|w|\theta_{R}^{1/2})^{2-\delta}\dd x\\
&\leq 2M_1^{\delta}||\nabla u_2||_{L^{r}}\Bigg(\int_{\Pi}|w|^{2}\theta_{R}\dd x \Bigg)^{1/r'}\\
&\leq 2M_1^{\delta}M_2 r \phi_{R}^{1/r'}.
\end{align*}\\
Thus, $\phi_R$ satisfies the differential inequality 

\begin{align*}
\dot{\phi_R}(t)&\leq a \phi_R(t)^{1/r'}+b,\quad t>0,\\
\phi_R(0)&=0,
\end{align*}\\
with the constants $a=2M_1^{\delta}M_2r$ and $b=C(R)$ by (5.6) and (5.7). Hence we have 

\begin{align*}
\int_{0}^{\phi_R(t)}\frac{\dd s}{as^{1/r'}+b}\leq t.  \tag{5.8}
\end{align*}\\

We prove that 

\begin{align*}
\overline{\lim}_{R\to\infty}\phi_R(t)<\infty\qquad \textrm{for each}\ t>0.   \tag{5.9}
\end{align*}\\
Suppose on the contrary that (5.9) were false for some $t_0>0$. Then, there exists a sequence $\{R_j\}$ such that $\lim_{j\to\infty}\phi_{R_j}(t_0)=\infty$. For an arbitrary $K>0$, we take a constant $N\geq 1$ such that $\phi_{R_j}(t_0)\geq K$ for $j\geq N$. It follows from (5.8) that 

\begin{align*}
\int_{0}^{K}\frac{\dd s}{as^{1/r'}+b}\leq t_0.  
\end{align*}\\
Since the constant $b=C(R_{j})$ converges to zero as $R_{j}\to\infty$, sending $j\to\infty$ yields $(r/a)K^{1/r}\leq t_0$. Since $K>0$ is arbitrary, this yields a contradiction. Thus (5.9) holds.

Since $|w|^{2}\theta_R$ monotonically converges to $|w|^{2}$ in $\Pi$, it follows from (5.9) that

\begin{align*}
\phi(t):=\int_{\Pi}|w|^{2}(x,t)\dd x=\lim_{R\to \infty}\int_{\Pi}|w|^{2}(x,t)\theta_{R}(x_3)\dd x<\infty.
\end{align*}\\
Sending $R\to\infty$ to (5.8) implies $(r/a)\phi^{1/r}(t)\leq t$. We thus obtain

\begin{align*}
\int_{\Pi}|w(x,t)|^{2}\dd x\leq M_1^{2}(2M_2 t)^{2/\delta}.
\end{align*}\\
Since the right-hand side converges to zero as $\delta\to0$ for $t\in [0,T]$ and $T=(4M_2)^{-1}$, we see that $w\equiv 0$ in $[0,T]$. Applying the same argument for $t\geq T$ implies $u_1\equiv u_2$ for all $t\geq 0$. The proof is now complete.
\end{proof}

\vspace{15pt}

\begin{rem}
By a similar cut-off function argument, uniqueness of weak solutions of the Euler equations with infinite energy is proved in \cite[Theorem 5.1.1]{Chemin} for the whole space under different assumptions from Theorem 1.2 (iii). See also \cite[Theorem 2]{Danchin}. We proved uniqueness of weak solutions in the infinite cylinder based on the Yudovich's estimate (Lemma 3.3).
\end{rem}

\vspace{15pt}

\section{Solutions with finite energy}

\vspace{15pt}

It remains to prove Theorem 1.3. The proof of the $L^{2}$-convergence (1.7) is simpler than that of uniqueness of weak solutions since solutions are with finite energy.

\subsection{Energy dissipation}

\begin{prop}
The assertion of Theorem 1.3 (i) holds.
\end{prop}

\vspace{5pt}

\begin{proof}
Let $u_0\in L^{p}_{\sigma}\cap L^{2}$ be an axisymmetric vector field without swirl such that $\omega^{\theta}_{0}/r\in L^{q}$ for $q\in [3/2,2]$ and $1/p=1/q-1/3$. By Theorem 1.1 and Proposition 2.2, there exists a unique global-in-time solution $u_{\nu}\in BC([0,\infty); L^{p}\cap L^{2})$ of (1.1) satisfying the energy equality (1.5). Since $\omega^{\theta}_{0}/r\in L^{q}$, applying Lemma 2.5 yields

\begin{align*}
\Big\|\frac{\omega^{\theta}}{r}\Big\|_{L^{2}}\leq \frac{C}{(\nu t)^{\frac{3}{2}(\frac{1}{q}-\frac{1}{2}) }}\Big\|\frac{\omega^{\theta}_{0}}{r}\Big\|_{L^{q}}\qquad t\geq 0,\ \nu>0.
\end{align*}\\
It follows that 

\begin{align*}
\nu\int_{0}^{T}||\nabla u_{\nu}||_{L^{2}}^{2}\dd t
=\nu\int_{0}^{T}||\omega^{\theta}||_{L^{2}}^{2}\dd t
\leq \nu\int_{0}^{T}\Big\|\frac{\omega^{\theta}}{r}\Big\|_{L^{2}}^{2}\dd t
\leq C\Big\|\frac{\omega^{\theta}_{0}}{r}\Big\|_{L^{q}}^{2}(\nu T)^{\frac{5}{2}-\frac{3}{q}}.
\end{align*}\\
Thus, (1.6) holds.
\end{proof}

\vspace{15pt}

\subsection{$L^{2}$-convergence}

We prove Theorem 1.3 (ii). We use uniform bounds for the viscosity $\nu>0$.

\vspace{15pt}

\begin{prop}
Let $u_0\in L^{p}_{\sigma}\cap L^{2}$ be an axisymmetric vector field without swirl such that $\omega^{\theta}_{0}/r\in L^{q}\cap L^{\infty}$ for $q\in [3/2,2]$ and $1/p=1/q-1/3$. Let $u_{\nu}\in BC([0,\infty);L^{p}\cap L^{2} )\cap C^{\infty}(\overline{\Pi}\times (0,\infty))$ be a solution of (1.1). Then, the estimates

\begin{align*}
||u_{\nu}||_{L^{\infty}}&\leq C\Big\|\frac{\omega^{\theta}_{0}}{r}\Big\|_{L^{r_0} },  \tag{6.1}\\
||\nabla u_{\nu}||_{L^{r}}&\leq C'r \Big\|\frac{\omega^{\theta}_{0}}{r}\Big\|_{L^{r}\cap L^{r_0} }, \tag{6.2}\\
||\nabla u_{\nu}||_{L^{2}}&\leq \Big\|\frac{\omega^{\theta}_{0}}{r}\Big\|_{L^{2} },  \tag{6.3}\\
\end{align*}
hold for $t>0$ and $3<r_0<r<\infty$ with some constants $C$ and $C'$, independent of $r$ and $\nu$.
\end{prop}

\vspace{5pt}

\begin{proof}
We take $r_0\in (3,\infty)$. It follows from (4.3) and (4.6) that 

\begin{align*}
||u_{\nu}||_{W^{1,r_0}}\leq C \Big\|\frac{\omega^{\theta}_{0}}{r}\Big\|_{L^{r_0}}\quad t\geq 0, \nu>0.
\end{align*}\\
By the Sobolev inequality, the estimate (6.1) follows. The estimates (6.2) and (6.3) follow from (2.5) and (3.15).
\end{proof}

\vspace{15pt}

Let $(u_{\nu}, p_{\nu})$ and $(u_{\mu}, p_{\mu})$ be two solutions of (1.1) for the same initial data $u_0$. We may assume that $\nu\geq \mu$. Then,  $w=u_{\nu}-u_{\mu}$ and $\pi=p_{\nu}-p_{\mu}$ satisfy

\begin{align*}
\partial_t w-\nu\Delta w-(\nu-\mu)\Delta u_{\mu}+u_{\nu}\cdot \nabla w+w\cdot \nabla u_{\mu}+\nabla \pi=0\quad \D\ w&=0\quad \textrm{in}\ \Pi\times (0,\infty),\\
\nabla \times w\times n=0,\ w\cdot n&=0\quad \textrm{on}\ \partial\Pi\times (0,\infty),\\
w&=0\quad \textrm{on}\ \Pi\times \{t=0\}.
\end{align*}\\
By multiplying $2w$ by the equation and integration by parts, we see that 

\begin{align*}
\frac{\dd}{\dd t}\int_{\Pi}|w|^{2}\dd x
+2\nu \int_{\Pi}|\nabla w|^{2}\dd x
+2(\nu-\mu)\int_{\Pi}\nabla u_{\mu}\cdot \nabla w\dd x
+2\int_{\Pi}(w\cdot \nabla u_{\mu})\cdot w\dd x=0.
\end{align*}\\
We set 

\begin{align*}
\phi_{\nu}(t)=\int_{\Pi}|w(x,t)|^{2}\dd x.
\end{align*}\\
We show that $K_{\nu}(T)=\sup_{0\leq t\leq T}\phi_{\nu}(t)$ converges to zero as $\nu \to 0$ for each $T>0$.

\vspace{15pt}

\begin{prop}
There exist constants $M_1-M_3$, independent of $\nu, \mu>0$ such that 

\begin{align*}
||w||_{L^{\infty}}&\leq M_1,\\
||\nabla u_{\mu}||_{L^{r}}&\leq M_2 r,\\
||\nabla u_{\mu}||_{L^{2}}+||\nabla w||_{L^{2}}&\leq M_3,\\
\end{align*}\\
hold for $t>0$ and $r>3$.
\end{prop}

\vspace{5pt}

\begin{proof}
The assertion follows from Proposition 6.2.
\end{proof}

\vspace{5pt}

\begin{proof}[Proof of Theorem 1.3 (ii)]
By Proposition 6.3, we estimate

\begin{align*}
\Big|2(\nu-\mu)\int_{\Pi}\nabla u_{\mu}\cdot \nabla w\dd x\Big| 
\leq 2\nu ||\nabla u_{\mu}||_{L^{2}}||\nabla w||_{L^{2}}
\leq 2\nu M_{3}^{2}.
\end{align*}\\
For an arbitrary $\delta \in (0, 2/3)$, we set $r=2/\delta$. Since $r'=2/(2-\delta)$, by a similar way as in the proof of Theorem 1.2 (iii), we estimate

\begin{align*}
\Big|2\int_{\Pi}(w\cdot \nabla u_{\mu})\cdot w\dd x\Big|
\leq 2M_1^{\delta}M_2 r\phi_{\nu}^{1/r'}.
\end{align*}\\
Thus, $\phi_{\nu}$ satisfies

\begin{align*}
&\dot{\phi_{\nu}}(t)\leq a\phi_{\nu}^{1/r'}(t)+b,\\
&\phi_{\nu}(0)=0,
\end{align*}\\
for $a=2M_1^{\delta}M_2r$ and $b=2\nu M_3^{2}$. We take an arbitrary $T>0$. We integrate the differential inequality between $(0,t)$ and take a supremum for $t\in [0,T]$ to estimate

\begin{align*}
\int_{0}^{K_{\nu}(T)}\frac{\dd s}{as^{1/r'}+b}\leq T.
\end{align*}\\
Since $b=b_{\nu}$ converges to zero as $\nu\to0$, by the same way as in the proof of Theorem 1.2 (iii), we see that the limit superior of $K_{\nu}(T)$ is finite for each $T>0$. We set 

\begin{align*}
K(T):=\overline{\lim}_{\nu\to0}K_{\nu}(T)<\infty.
\end{align*}\\
Sending $\nu\to0$ to the above inequality implies $(r/a)K^{1/r'}(T)\leq T$. Since $a=2M_1^{\delta}M_2r$ and $r=2/\delta$, it follows that

\begin{align*}
K(T)\leq M_1^{2}(2M_2T)^{2/\delta}.
\end{align*}\\
Since the right-hand side converges to zero as $\delta\to0$ for $T\leq T_0$ and $T_0=(4M_2)^{-1}$, we see that $K(T)\equiv 0$. Thus the convergence (1.7) holds. By replacing the initial time and applying the same argument for $T\geq T_0$, we are able to show the convergence (1.7) for an arbitrary $T>0$. The proof is now complete.
\end{proof}

\vspace{15pt}

\begin{rems}
\noindent 
(i) (Convergence in Sobolev space) 
We constructed global weak solutions of the Euler equations for axisymmetric data without swirl by a vanishing viscosity method. As explained in the introduction, the condition $\omega^{\theta}_{0}/r\in L^{q}$ for $q\in [3/2,3)$ in Theorem 1.2 (i) is satisfied if $u_0\in W^{2,q}$. This condition is weaker than $u_0\in W^{2,q}$ for $q\in (3,\infty)$, required for the local well-posedness of the Euler equations. 

Our approach is based on the a priori estimate (1.10) which is a special property of axisymmetric solutions and is not available at the broad level. On the other hand, there is an another approach to study vanishing viscosity limits when the Euler equation is locally well-posed. When $\Pi=\mathbb{R}^{3}$, unique local-in-time solutions of the Euler equations are constructed in \cite{Swann}, \cite{Kato72}, \cite{Kato75} by sending $\nu\to0$ to local-in-time solutions of the Navier-Stokes equations. See also \cite{Constantin86}. In particular, for a local-in-time solution $u\in C([0,T]; H^{s})$ of the Euler equations and $u_0\in H^{s}$, $s>5/2$, the convergence 

\begin{align*}
u_{\nu} \to u\quad \textrm{in}\ L^{\infty}(0,T; H^{s}),
\end{align*}\\
is known to hold \cite{Masmoudi07}. The case with boundary is a difficult question related to analysis of boundary layer. See \cite{Constantin07} for a survey. However, convergence results are known subject to the Neumann boundary condition (1.1). See \cite{XiaoXin}, \cite{Veiga10}, \cite{Veiga11} for the case with flat boundaries and \cite{BS12}, \cite{BS14} for curved boundaries.

\noindent
(ii) (Navier boundary condition) The Neumann boundary condition in (1.1) may be viewed as a special case of the Navier boundary condition,

\begin{align*}
(D(u)n+\alpha u)_{\textrm{tan}}=0,\quad u\cdot n=0\quad \textrm{on}\ \partial\Pi,  \tag{6.4}
\end{align*}\\
where $D(u)=(\nabla u+\nabla^{T}u)/2$ is the deformation tensor and $f_{\textrm{tan}}=f-n(f\cdot n)$ for a vector field $f$. Indeed, for the two-dimensional case, the Neumann boundary condition is reduced to the free condition $\omega=0$ and $u\cdot n=0$ on $\partial\Pi$. The free condition is a special case of (6.4), which is written as $\omega+2(\alpha-\kappa)u\cdot n^{\perp}=0$ and $u\cdot n=0$ on $\partial\Pi$, with the curvature $\kappa(x)$ and $n^{\perp}=(-n^{2},n^{1})$. For a two-dimensional bounded domain, vanishing viscosity limits subject to (6.4) are studied in \cite{Robert98}, \cite{Planas05}, \cite{Ke06}. For the three-dimensional case, it is shown in \cite{IP06} that a Leray-Hopf weak solution $u_{\nu}$ subject to (6.4) converges to the local-in-time solution $u\in C([0,T]; H^{s})$ of the Euler equations for $u_0\in H^{3}$ in the sense that

\begin{align*}
u_{\nu}\to u\quad \textrm{in}\ L^{\infty}(0,T; L^{2}).
\end{align*}\\
For the Dirichlet boundary condition, the same convergence seems unknown. See \cite{IS11} for a boundary layer expansion subject to (6.4) and \cite{MR12} for a stronger convergence result.
\end{rems}

\vspace{30pt}

\appendix

\section{Decay estimates of vorticity}

\vspace{15pt}

We prove the decay estimate (2.9) (Lemma 2.5). It suffices to show:

\vspace{15pt}

\begin{lem}
There exists a constant $C$ such that the estimate (2.9) holds for $r=2^{m}q$, $q\in [1,\infty)$ and non-negative integers $m\geq 0$.
\end{lem}

\vspace{5pt}

\begin{proof}[Proof of Lemma 2.5]
We apply Lemma A.1. Since 

\begin{align*}
\lim_{r\to\infty}\Big\|\frac{\omega^{\theta}}{r}\Big\|_{L^{r}}(t)=\Big\|\frac{\omega^{\theta}}{r}\Big\|_{L^{\infty}}(t),
\end{align*}\\
sending $m\to\infty$ implies (2.9) for $r=\infty$ and $q\in [1,\infty)$. Since (2.9) holds for $r=q\in [1,\infty]$, we obtain (2.9) for all $1\leq q\leq r\leq \infty$ by the H\"older inequality.
\end{proof}

\vspace{15pt}

Let $\psi_{\varepsilon}(s)$ be a non-negative convex function in the proof of Proposition 2.4. We prove the estimate (2.9) for $\psi_{\varepsilon}(\Omega)$ and $\Omega=\omega^{\theta}/r$. The assertion of Lemma A.1 follows by sending $\varepsilon\to0$.

\vspace{15pt}

\begin{prop}
There exists a constant $C$ such that the estimate 

\begin{align*}
\big\|\psi_{\varepsilon}(\Omega)\big\|_{L^{r}(\Pi)}\leq \frac{C}{(\nu t)^{\frac{3}{2}(\frac{1}{q}-\frac{1}{r}) }}\big\|\psi_{\varepsilon}(\Omega_0)\big\|_{L^{q}(\Pi)}\qquad t>0,\ \nu>0,  \tag{A.1}
\end{align*}\\
holds for all $\varepsilon>0$, $r=2^{m}q$ and $m\geq 0$.
\end{prop}

\vspace{15pt}

We consider differential inequalities for $L^{r}$-norms of $\psi_{\varepsilon}(\Omega)$.

\vspace{15pt}

\begin{prop}
The function 

\begin{align*}
\phi_{r}(t)=\int_{\Pi}\psi_{\varepsilon}^{r}(\Omega)\dd x
\end{align*}\\
satisfies 

\begin{align*}
\dot{\phi}_{r}(t)\leq -\kappa \nu \Big(1-\frac{1}{r}\Big)\frac{\phi_{r}^{5/3}}{\phi_{r/2}^{4/3}}\qquad t>0    \tag{A.2}
\end{align*}\\
with an absolute constant $\kappa$, independent of $r$ and $\nu$.
\end{prop}

\vspace{15pt}

\begin{proof}
We apply the interpolation inequality 

\begin{align*}
||\varphi||_{L^{2}}\leq C_0||\varphi||_{L^{1}}^{2/5}||\nabla \varphi||_{L^{2}}^{3/5}
\end{align*}\\
for $\varphi\in H^{1}_{0}$ with some absolute constant $C_{0}$. Since $\psi_{\varepsilon}(\Omega)$ satisfies 

\begin{align*}
\frac{\dd}{\dd t}\int_{\Pi}\psi^{r}_{\varepsilon}(\Omega)\dd x
+4\nu \Big(1-\frac{1}{r}\Big) \int_{\Pi}\Big|\nabla \psi_{\varepsilon}(\Omega)^{\frac{r}{2}}\Big|^{2}\dd x\leq 0,
\end{align*}\\
by (2.7), applying the interpolation inequality for $\varphi=\psi_{\varepsilon}^{r/2}$ yields 

\begin{align*}
\int_{\Pi}\Big|\nabla \psi_{\varepsilon}(\Omega)^{\frac{r}{2}}\Big|^{2}\dd x\geq \frac{1}{C_0^{10/3}}
\frac{\Big(\int_{\Pi}\psi_{\varepsilon}^{r}\dd x\Big)^{5/3}}{\Big(\int_{\Pi}\psi_{\varepsilon}^{r/2}\dd x\Big)^{4/3}}=\frac{\phi_{r}^{5/3}}{C_0^{10/3}\phi_{r/2}^{4/3}}.
\end{align*}\\
The differential inequality (A.2) follows from the above two inequalities with $\kappa=4C_{0}^{-10/3}$.
\end{proof}

\vspace{15pt}

\begin{proof}[Proof of Proposition A.2]
We set $\lambda=||\psi_{\varepsilon}(\Omega_0)||_{L^{q}}$. The estimate (A.1) is written as

\begin{align*}
\phi_{r}^{1/r}(t)\leq \frac{C}{(\nu t)^{\frac{3}{2}(\frac{1}{q}-\frac{1}{r}) }}\lambda\qquad t>0, \tag{A.3}
\end{align*}\\
for $r=2^{m}q$ and $m\geq 0$. We prove (A.3) by induction for $m\geq 0$. For $m=0$, the estimate (A.3) holds with $C=1$ by Lemma 2.3.

Suppose that (A.3) holds for $m=k$ with some constant $C=C_k$. We set $s=2r$ for $r=2^{m}q$. By the assumption of our induction, we see that 

\begin{align*}
\frac{1}{\phi_r(t)}\geq \frac{(\nu t)^{\frac{3}{2}(2^{k}-1) }}{C_{k}^{r}\lambda^{r}}.
\end{align*}\\
It follows from (A.2) that 

\begin{align*}
\phi_{s}^{-5/3}\dot{\phi}_{s}
&\leq -\kappa \nu \Big(1-\frac{1}{s}\Big)\phi_{r}^{-4/3}\\
&\leq -\kappa \nu\Big(1-\frac{1}{s}\Big) \frac{(\nu t)^{2^{k+1}-2 }}{C_{k}^{4r/3}\lambda^{4r/3} }.
\end{align*}\\
We integrate the both sides between $[t_1,t]$ and estimate 

\begin{align*}
\frac{3}{2}\Bigg(\frac{1}{\phi^{2/3}(t)}-\frac{1}{\phi^{2/3}(t_1)}  \Bigg)
\geq \kappa \Big(1-\frac{1}{s}\Big) \frac{(\nu t)^{2^{k+1}-1 }-(\nu t_1)^{2^{k+1}-1 }}{C_{k}^{4r/3}\lambda^{4r/3} }\Bigg(\frac{1}{2^{k+1}-1}\Bigg).
\end{align*}\\
Since the left-hand side is smaller than $3/2\phi_{s}^{-2/3}$, sending $t_1\to0$ yields 

\begin{align*}
\phi_s^{2/3}(t)\leq \frac{(C_k\lambda)^{4r/3} }{(\nu t)^{2^{k+1}-1} }\frac{ 2^{k+2}}{\kappa (1-1/s)}.
\end{align*}\\
Since 

\begin{align*}
\frac{3}{2s}\frac{4r}{3}=1,\qquad
\frac{3}{2s}(2^{k+1}-1)=\frac{3}{2}\Big(\frac{1}{q}-\frac{1}{s}\Big),
\end{align*}\\
and $1-1/s\geq 1/2$, it follows that 

\begin{align*}
\phi_{s}^{1/s}(t)\leq \frac{C_k\lambda }{(\nu t)^{\frac{3}{2}(\frac{1}{q}-\frac{1}{s})} }\Bigg(\frac{ 2^{k+3}}{\kappa}\Bigg)^{\frac{3}{2^{k+2}q}}.
\end{align*}\\
We proved that (A.3) holds for $m=k+1$ with the constant $C_{k+1}=a_kC_k$ for $a_{k}=b^{\frac{1}{2^{k+2}}}d^{\frac{k+3}{2^{k+2}}}$ and

\begin{align*}
b=\kappa^{-3/q},\qquad d=2^{3/q}.
\end{align*}\\
Thus (A.3) holds for all $m\geq 0$. Since

\begin{align*}
C_{k+1}=a_kC_k
=\prod_{j=1}^{k}a_{j}
=b^{\sum_{j=1}^{k}2^{-j-2}}d^{\sum_{j=1}^{k}(j+3)2^{-j-2}},
\end{align*}\\
and the right-hand side converges as $k\to \infty$, we are able to take a uniform constant $C$ in (A.3) for all $m\geq 0$. The proof is now complete.
\end{proof}

\vspace{15pt}

\section*{Acknowledgements}
The author is partially supported by JSPS through the Grant-in-aid for Young Scientist (B) 17K14217, Scientific Research (B) 17H02853 and Osaka City University Strategic Research Grant 2018 for young researchers.

\vspace{10pt}


\vspace{10pt}

\bibliographystyle{plain}
\bibliography{ref}

\end{document}